\documentclass[11pt, letterpaper]{amsart}
\usepackage{amstext,amsthm,amsmath,amsfonts,amssymb,amscd,latexsym,txfonts,stmaryrd}
\numberwithin{equation}{section}
\usepackage{color}
\usepackage{epsfig}
\usepackage{cite}
\usepackage[colorlinks=true, pdfstartview=FitV, linkcolor=blue, citecolor=blue, urlcolor=blue]{hyperref}
\usepackage[nameinlink]{cleveref}
\usepackage[all,cmtip]{xy}
\usepackage{youngtab}
\usepackage{enumerate}

\newcommand{\dsp}{\displaystyle}

\newcommand{\C}{{\mathbb C}}
\newcommand{\R}{{\mathbb R}}
\newcommand{\Z}{{\mathbb Z}}

\newcommand{\N}{{\mathbb N}}

\newcommand{\GL}{\operatorname{GL}}

\newcommand{\Ann}{\operatorname{Ann}}
\newcommand{\flo}[1]{\left\lfloor\frac{#1}{2}\right\rfloor}

\newcommand{\Hess}{\operatorname{Hess}}
\newcommand{\sgn}{\operatorname{sgn}}

\newcommand{\supp}{\operatorname{supp}}
\newcommand{\Lor}{\overset{\circ}{L}}

\usepackage{xcolor}

\theoremstyle{theorem}
\newtheorem{theorem}{Theorem}[section]
\newtheorem{proposition}[theorem]{Proposition}
\newtheorem{lemma}[theorem]{Lemma}
\newtheorem{corollary}[theorem]{Corollary}
\newtheorem{problem}[theorem]{Problem}

\newtheorem{introthm}{Theorem}
\newtheorem{introcor}{Corollary}

\newtheorem{fact}[theorem]{Fact}

\newtheorem*{question*}{Question}
\newtheorem*{answer*}{Answer}
\newtheorem*{remarks*}{Remarks}
\newtheorem*{claim*}{Claim}
\newtheorem*{remark*}{Remark}
\newtheorem*{proposition*}{Proposition}
\newtheorem*{lemma*}{Lemma}
\newtheorem*{fact*}{Fact}

\theoremstyle{definition}
\newtheorem{definition}[theorem]{Definition}

\newtheorem{remark}[theorem]{Remark}
\newtheorem{example}[theorem]{Example}

\begin{document}
	\title[Higher Lorentzian Polynomials]{Higher Lorentzian Polynomials, Higher Hessians, and the Hodge-Riemann Relations for graded oriented Artinian Gorenstein Algebras in Codimension Two}
	\author[MMMSW]{Pedro Macias Marques, Chris McDaniel, Alexandra Seceleanu, Junzo Watanabe}
	
	\address{Departamento de Matem\'{a}tica \\
		Escola de Ci\^{e}ncias e Tecnologia, Centro de Investiga\c{c}\~{a}o em Matem\'{a}tica e Aplica\c{c}\~{o}es\\
		Instituto de Investiga\c{c}\~{a}o e Forma\c{c}\~{a}o Avan\c{c}ada\\
		Universidade de \'{E}vora\\
		Rua Rom\~{a}o Ramalho, 59, P--7000--671 \'{E}vora, Portugal}
	\email{pmm@uevora.pt}
	
	\address{Department of Mathematics\\
		Endicott College\\
		376 Hale St Beverly, MA 01915, USA.}
	\email{cmcdanie@endicott.edu}

	\address{Department of Mathematics\\
		University of Nebraska-Lincoln\\
		Lincoln NE 68588, USA.}
	\email{aseceleanu@unl.edu}
	
	\address{Department of Mathematics\
		Tokai University\\
		Hiratsuka, Kanagawa 259-1292 Japan}
	\email{watanabe.junzo@tokai-u.jp}

	
	\maketitle
	\begin{abstract}
		A (standard graded) oriented Artinian Gorenstein algebra over the real numbers is uniquely determined by a real homogeneous polynomial called its Macaulay dual generator.  We study the mixed Hodge-Riemann relations on oriented Artinian Gorenstein algebras for which we give a signature criterion on the higher mixed Hessian matrices of its Macaulay dual generator.  Inspired by recent work of Br\"and\'en and Huh, we introduce a class of homogeneous polynomials in two variables called $i$-Lorentzian polynomials, and show that these are exactly the Macaulay dual generators of oriented Artinian Gorenstein algebras in codimension two satisfying mixed Hodge-Riemann relations up to degree $i$ on the positive orthant of linear forms.  We further show that the set of $i$-Lorentzian polynomials of degree $d$ are in one-to-one correspondence with the set of totally nonnegative Toeplitz matrices of size depending on $i$ and $d$.  A corollary is that all normally stable polynomials, i.e. polynomials whose normalized coefficients form a PF sequence, are $i$-Lorentzian.  Another corollary is an analogue of Whitney's theorem for Toeplitz matrices, which appears to be new:  the closure of the set of totally positive Toeplitz matrices, in the Euclidean space of all real matrices of a given size, is equal to the set of totally nonnegative Toeplitz matrices.  
	\end{abstract}

\section{Introduction}
	
	
In their recent paper \cite{BH}, Br\"and\'en-Huh defined a remarkable class of homogeneous polynomials in $n$ variables that they called Lorentzian polynomials.  Lorentzian polynomials provide a direct link between the Hodge-Riemann relations in degree one, an algebraic condition motivated by K\"ahler geometry, and the combinatorial condition of log concavity.  In this paper, we focus on the $n=2$ variables case, and define higher Lorentzian polynomials which connect Hodge-Riemann relations in higher degrees with certain higher log concavity conditions encapsulated in the total positivity or nonnegativity of certain Toeplitz matrices.

A graded oriented Artinian Gorenstein (AG) algebra $A$ satisfies the strong Lefschetz property in degree $i$ if there is some linear form $\ell\in A_1$ such that the multiplication maps $\times\ell^{d-2j}\colon A_j\rightarrow A_{d-j}$ are vector space isomorphisms for each $0\leq j\leq i$, and it satisfies the Hodge-Riemann relations in degree $i$ if those same multiplication maps additionally satisfy a certain alternating signature condition (\Cref{def:HRR}).  More generally, $A$ satisfies the mixed strong Lefschetz property or the mixed Hodge-Riemann relations on a subset of linear forms $U\subset A_1$ if we replace multiplication by powers of a single linear form with multiplication by arbitrary products of linear forms from $U$ (\Cref{def:MHRR}). 
The strong Lefschetz property and the Hodge-Riemann relations, and their mixed analogues, are well known properties of the cohomology ring of a compact complex K\"ahler manifold, but they are also known to hold in various other non-geometric settings, where they have been used to establish various combinatorial inequalities:  unimodality of the $h$-vectors of simple polytopes \cite{S,McM,T}, positivity of Kazhdan-Lusztig coefficients \cite{EW}, and log concavity of the characteristic polynomial of a matroid \cite{AHK}, to name a few.  Whereas the strong Lefschetz property implies certain linear inequalities on the Hilbert function of an algebra, the Hodge-Riemann relations imply certain higher degree polynomial inequalities on the coefficients of the Macaulay dual generator. 

Every graded oriented AG algebra of socle degree $d$ is uniquely determined by a homogeneous polynomial of degree $d$.  The algebra is denoted by $A_F$ and $F$ is called its Macaulay dual generator, named after work of F.S. Macaulay \cite{Macaulay} in the early 20${^{\text{th}}}$ century; this polynomial is also known by some authors as the volume polynomial of $A_F$, e.g. \cite{AHK,BH,H,T}.  The matrices for the Lefschetz multiplication maps, with respect to a certain basis $\mathcal{E}$ of the algebra $A_F$, can be identified with the higher Hessian matrices of $F$ (\Cref{fact:HessLef}), a fact discovered by the fourth author \cite{W}; see also \cite{MW}.  We show that the Hodge-Riemann relations are equivalent to certain signature conditions on these higher Hessian matrices (\Cref{fact:HessHRR}), and the mixed Hodge-Riemann relations are equivalent to those same signature conditions on certain polarized versions of the higher Hessian matrices that we call mixed Hessians (\Cref{fact:MHEssHRR}).  

It is well known that AG algebras in codimension two (corresponding to Macaulay dual generators in $n=2$ variables) are special among AG algebras.  For example, AG algebras in codimension two are always complete intersections \cite{IK} and they always satisfy the strong Lefschetz property \cite{I}; however, as we shall see, they do not always satisfy the Hodge-Riemann relations (\Cref{ex:SLnotHRR}).  Furthermore, in the codimension two case, we show that the aforementioned Hodge-Riemann signature conditions amount to certain alternating sign conditions on the determinants of these higher Hessians and their mixed analogues (\Cref{fact:TwoMHessHRR}).  Moreover in this case, each (mixed) Hessian matrix is Hankel, and can be transformed into a Toeplitz matrix by multiplying by an appropriate permutation matrix, whereupon the alternating sign condition is transformed into a positivity condition (\Cref{lem:TwoHess}).  After some generalities, we shall focus exclusively on the codimension two case, where we characterize oriented AG algebras satisfying mixed Hodge-Riemann relations in terms of their Macaulay dual generators and the total positivity of their associated Toeplitz matrices.  Our main results are stated below.

Given integers $i,d$ satisfying $0\leq i\leq \flo{d}$, there is a linear isomorphism $\phi^i_d\colon \R[X,Y]_d\rightarrow\mathcal{T}(i+1,d-i+1)$ between the set of bivariate homogeneous polynomials of degree $d$ and the set of real $(i+1)\times (d-i+1)$ Toeplitz matrices.  We define the set of strictly $i$-Lorentzian polynomials as the set of polynomials $F$ in which all (consecutive) minor determinants of the Toeplitz matrix $\phi^i_d(F)$ are positive, and, following Br\"and\'en and Huh \cite{BH}, we define the set of $i$-Lorentzian polynomials to be their closure in the Euclidean space $\R[X,Y]_d\cong \R^{d+1}$ (\Cref{def:Lorentzian}).  Recall that a matrix is totally positive, respectively nonnegative, if all of its minor determinants are positive, respectively nonnegative. 
\begin{introthm}
	\label{introthm:strictLor}
	The following are equivalent:
	\begin{enumerate}
		\item $F$ is strictly $i$-Lorentzian
		\item the Toeplitz matrix $\phi^i_d(F)$ is totally positive
		\item the oriented AG algebra $A_F$ satisfies mixed Hodge-Riemann relations in degree $i$ on the standard closed convex cone  
		$$\overline{U}=\left\{ax+by \ | \ (a,b)\in\R^2_{\geq 0}\setminus\{(0,0)\}\right\}$$
	\end{enumerate}
\end{introthm}

\begin{introthm}
	\label{introthm:Lor}
	The following are equivalent:
	\begin{enumerate}
		\item $F$ is $i$-Lorentzian
		\item the Toeplitz matrix $\phi^i_d(F)$ is totally nonnegative
		\item the oriented AG algebra $A_F$ satisfies mixed Hodge-Riemann relations in degree $i$ on the standard open convex cone 
		$$U=\left\{ax+by \ | \ (a,b)\in\R^2_{>0}\right\}$$
	\end{enumerate}
\end{introthm}
In the case $i=1$, $\phi^1_d(F)$ is a two-rowed Toeplitz matrix whose total nonnegativity is equivalent to the coefficient sequence of $F$ being nonnegative, and ultra log concave with no internal zeros.  In particular our $1$-Lorentzian polynomials agree with the Lorentzian polynomials defined by Br\"and\'en and Huh \cite{BH} in $n=2$ variables (\Cref{fact:BHLor}).  It would be interesting to characterize $i$-Lorentzian polynomials in $n>2$ variables, perhaps taking \Cref{introthm:Lor}(3) as a definition (\Cref{prob:2}).  As a corollary to \Cref{introthm:Lor}, we obtain the following purely linear algebraic result, which appears to be new:
\begin{introcor}
	\label{introthm:Toep}
	The closure of the set of totally positive Toeplitz matrices is equal to the set of totally nonnegative Toeplitz matrices. 
\end{introcor} 

\Cref{introthm:Toep} is an analogue of a result of A. Whitney \cite{Wh}; see also \cite[Theorem 2.7]{Ando}. 
The idea of the proof of \Cref{introthm:Lor}, and hence also \Cref{introthm:Toep}, is to realize the permuted higher (mixed) Hessian matrix of $F$ as a certain finite minor of the matrix product of two infinite Toeplitz matrices, one of which turns out to be the totally nonnegative weighted path matrix of a certain edge-weighted acyclic planar graph, a subgraph of the NE lattice paths in the plane (\Cref{thm:HessMHess}).  For further connections between weighted path matrices and total nonnegativity, see the survey article by Fomin-Zelevinsky \cite{FZ}.

Closely related to the notion of Lorentzian polynomials is that of stable polynomials; in two variables, these are real rooted polynomials with nonnegative coefficients.  Br\"and\'en and Huh \cite{BH} showed that stable polynomials are always $1$-Lorentzian, however, as we shall see, they need not always be $i$-Lorentzian for $i>1$ (\Cref{ex:NSLor}).  We define a subclass of stable polynomials called normally stable polynomials, essentially obtained from stable polynomials by multiplying their coefficients by the binomial coefficients (\Cref{def:stable}).

A beautiful characterization of stable polynomials in two variables was obtained by Aissen-Schoenberg-Whitney  \cite{ASW} and Edrei \cite{E} in the 1950s (\Cref{fact:ASWE}).  From their characterization it follows directly that a homogeneous bivariate polynomial is normally stable if and only if its \emph{infinite} Toeplitz matrix $\phi_d(F)$ is totally nonnegative, of which the finite Toeplitz matrices $\phi^i_d(F)$ above are minors.  This fact, taken together with our characterization of bivariate $i$-Lorentzian polynomials \Cref{introthm:Lor}, gives the following result:
\begin{introcor}
	\label{introthm:Norm}
	Every normally stable polynomial is $i$-Lorentzian for all $i\geq 0$.
\end{introcor}  
\Cref{introthm:Norm} gives an easy way to build oriented AG algebras of codimension two that satisfy the Hodge-Riemann relations in any degree (\Cref{ex:NSLor}).   

This paper is organized as follows.  In \Cref{sec:Pre}, we discuss (mixed) Hodge-Riemann relations on oriented AG algebras.  In \Cref{sec:MDHess}, we discuss Macaulay duality, higher (mixed) Hessians, and state the (mixed) Hodge-Riemann relations signature conditions on them.  In \Cref{sec:HLor}, we define strictly $i$-Lorentzian, $i$-Lorentzian, and normally stable polynomials, and prove \Cref{introthm:strictLor} (\Cref{thm:Lorphi} and \Cref{thm:LorMHRR}), \Cref{introthm:Lor} (\Cref{thm:Lorentzian}), \Cref{introthm:Toep} (\Cref{cor:Toep}), and \Cref{introthm:Norm} (\Cref{thm:NStable}).  In \Cref{sec:Conclude} we discuss some open problems.  In \Cref{sec:App} we give a brief review of some results of Br\"and\'en-Huh \cite{BH}.  In this paper, all algebras are assumed to be standard graded over $\R$, the field of real numbers, unless otherwise stated.

	\section{Preliminaries}
	\label{sec:Pre}
	\subsection{Oriented AG Algebras and Hodge-Riemann Relations}
	An oriented Artinian Gorenstein (AG) $\R$-algebra is a graded AG algebra $A=\bigoplus_{i=0}^dA_i$ with $A_0=\R$ together with a fixed linear isomorphism
	$$\int_A\colon A_d\rightarrow \R$$
	such that the \emph{intersection pairing} in each degree $i$ defined by multiplication in $A$ is nondegenerate:
	$$\xymatrixrowsep{.5pc}\xymatrix{A_i\times A_{d-i}\ar[r] & \R\\
		(\alpha,\beta)_i\ar@{=}[r] & \int_A\alpha\beta.\\}$$
	
	The \emph{Hilbert function} of $A$ is the finite integer sequence $H(A)=(h_0(A),h_1(A),\ldots,h_d(A))$ where $h_i(A)=\dim_\R(A_i)$; as a consequence of the nondegenerate intersection pairing, the Hilbert function of an AG algebra is always symmetric, i.e. $h_i(A)=h_{d-i}(A)$ for all $0\leq i\leq d$. 
	
	To a linear form $\ell\in A_1$, and an integer $0\leq i\leq\flo{d}$, we define the $i^{th}$ primitive subspace with respect to $\ell$
	\[P_{i,\ell}=\ker\left(\times \ell^{d-2i+1}\colon A_i\rightarrow A_{d-i+1}\right)\]
	and the $i^{th}$ Lefschetz form with respect to $\ell$
	$$\xymatrixrowsep{.5pc}\xymatrix{A_i\times A_i\ar[r] & \R\\ (\alpha,\beta)_i^\ell\ar@{=}[r] & (-1)^i\int_A\ell^{d-2i}\alpha\beta\\}.$$
	\begin{definition}
		\label{def:HRR}
		The pair $(A,\ell)$ is strong Lefschetz in degree $i$ (SL$_i$) if for each $0\leq j\leq i$, the $j^{th}$ Lefschetz form with respect to $\ell$ is non-degenerate, or equivalently, the $j^{th}$ Lefschetz multiplication map $\times\ell^{d-2j}\colon A_j\rightarrow A_{d-j}$ is an isomorphism.  It satisfies the strong Lefschetz property (SLP) if it satisfies SL$_{\flo{d}}$.
		
		The pair $\left(\left(A,\int_A\right),\ell\right)$ satisfies the (ordinary) Hodge-Riemann relations in degree $i$ (HRR$_i$)
		if for every $0\leq j\leq i$, the $j^{th}$ Lefschetz form with respect to $\ell$ is positive definite on the $j^{th}$ primitive subspace with respect to $\ell$, i.e. 
		$$(\alpha,\alpha)^\ell_j>0, \ \begin{array}{l} \forall 0\leq j\leq i\\ \forall 0\neq \alpha\in P_{j,\ell}\\ \end{array}$$
		The pair satisfies the Hodge-Riemann property (HRP) if it satisfies HRR$_{\flo{d}}$.
		For a subset $U\subset A_1$, we say that $\left(A,\int_A\right)$ satisfies HRR$_i$ on $U$ if $\left(\left(A,\int_A\right),\ell\right)$ satisfies HRR$_i$  for all $\ell\in U$.
	\end{definition}
\begin{remark}
	\label{rem:light}
	We remark that our notation for SL$_i$ may differ slightly from other notions found elsewhere in the literature; we have chosen it to be consistent with our notation for HRR$_i$, which includes all degrees \emph{up to} degree $i$.  Also we include the adjective \emph{ordinary} to distinguish it from the \emph{mixed} version that we introduce next; we shall drop this adjective if it is clear from the context.
\end{remark}

	\begin{lemma}
			\label{lem:HRRSL}
			If $\left(\left(A,\int_A\right),\ell\right)$ satisfies HRR$_i$, then $(A,\ell)$ also satisfies SL$_i$.
		\end{lemma}
	\begin{proof}
		Assume that $(A,\ell)$ does not satisfy SL$_i$.  Then for some $0\leq j\leq i$, there exists nonzero $\alpha\in \ker(\times\ell^{d-2j}\colon A_j\rightarrow A_{d-j})$.  Then $0\neq \alpha\in P_{j,\ell}$, and we have 
		$$\int_A\ell^{d-2j}\alpha^2=0$$
		which means that $\left(\left(A,\int_A\right),\ell\right)$ does not satisfy HRR$_i$. 
	\end{proof}
The converse of \Cref{lem:HRRSL} is not true, as the following example shows.
	\begin{example}\label{ex:SLnotHRR}
	Take $A=\R[x,y]/(x^2-y^2,xy)$ with $\int_A x^2=1$.  Then $\ell(a,b)=ax+by\in A_1$ satisfies SL$_1$ (hence also SLP) for every $(0,0)\neq (a,b)\in\R^2$, but $\left(\left(A,\int_A\right),\ell(a,b)\right)$ does not satisfy HRR$_1$ for any $(a,b)\in\R^2$ at all.
\end{example}

The following is sometimes referred to as the primitive decomposition of $A$ (with respect to $\ell$), and it is a direct consequence of SL$_i$.
	\begin{lemma}
		\label{lem:primitiveD}
		If $(A,\ell)$ satisfies SL$_i$, then for each $0\leq j\leq i$, $\dim_\R(P_{j,\ell})=h_j-h_{j-1}$ and there is an orthogonal decomposition with respect to the $j^{th}$ Lefschetz form
		$$A_j=P_{j,\ell}\oplus \ell(A_{j-1}).$$
	\end{lemma}
	\begin{proof}
	Fix $0\leq j\leq i$, fix $\alpha\in A_j$ and consider the element $y=\ell^{d-2j+1}\alpha\in A_{d-j+1}$.  If $y=0$, then $\alpha\in P_{j,\ell}$.  Otherwise, we may assume $j>0$ and $y\neq 0$ and by SL$_{j-1}$ there exists $\beta\in A_{j-1}$ such that $\ell^{d-2j+2}\beta=\ell^{d-2j+1}\alpha$, which implies that $\alpha-\ell\beta\in P_{j,\ell}$.  This shows that $A_j=P_{j,\ell}+\ell(A_{j-1})$.  It remains to see that the sum is direct and that the decomposition is orthogonal with respect to the $j^{th}$ Lefschetz pairing.  Suppose that $\alpha\in P_{j,\ell}\cap \ell(A_{j-1})$.  Then $\alpha=\ell\beta$ and $\ell^{d-2j+1}\alpha=\ell^{d-2j+2}\beta=0$, and thus by SL$_{j-1}$ we must have $\beta=0$ and hence $\alpha=0$.  If $\alpha\in P_{j,\ell}$ and $\ell\beta\in \ell(A_{j-1})$ then $$(\alpha,\ell\beta)_j^\ell=(-1)^j\int_A\ell^{d-2j}\alpha\ell\beta=(-1)^j\int_A\ell^{d-2j+1}\alpha\beta=0,$$
		as desired.  Finally it follows from SL$_i$ that the multiplication maps $\times\ell\colon A_{j-1}\rightarrow A_j$ are injective for all $0\leq j\leq i$, hence $\dim_\R(P_{j,\ell})=\dim_\R(A_j)-\dim_\R(\ell(A_{j-1}))=h_j-h_{j-1}$.     
	\end{proof}

	The following result will be useful for inductive arguments; it might be considered an algebraic analogue of the Lefschetz hyperplane theorem in algebraic geometry.
	
	\begin{lemma}
		\label{lem:Key}
		Let $\left(A,\int_A\right)$ be an oriented AG algebra, $\ell\in A_1$ a non-zero linear form, and consider the quotient $B=A/(0:\ell)$ by the colon ideal $(0:\ell)=\left\{a\in A \ | \ \ell a=0\right\}$ with $\pi\colon A\rightarrow B$ the natural surjection.  Then $B$ is an oriented AG algebra with orientation $\int_B\pi(\alpha) = \int_A\ell\alpha$.  Moreover, if the pair $\left(\left(A,\int_A\right)\ell\right)$ satisfies HRR$_i$ or HRP, then so does the pair $\left(\left(B,\int_B\right),\pi(\ell)\right)$ and for each $0\leq j\leq \min\Bigl\{i,\flo{d-1}\Bigr\}$ we have 
		$h_j(A)=h_j(B)$. 
	\end{lemma}
\begin{proof}
	To see that $B$ is Gorenstein, it suffices to see that the socle $$\operatorname{soc}(B)=\{b\in B \ | \ xb=0, \ \forall x\in\mathfrak{m}_B\}$$ 
	is one dimensional.  If $0\neq \beta\in \operatorname{soc}(B)$, then since $\pi(\mathfrak{m}_A)=\mathfrak{m}_B$, it follows that any lift $\hat{\beta}\in A$ must satisfy $0\neq \ell\hat{\beta}\in\operatorname{soc}(A)$.  It follows that $\operatorname{soc}(B)\cong \ell\left(\pi^{-1}(\operatorname{soc}(B))\right)\subseteq \operatorname{soc}(A)$.  Since $\operatorname{soc}(A)$ is one dimensional it follows that $\operatorname{soc}(B)$ is too, 
	and the map $\int_B=\int_A\ell(\cdot)\colon \operatorname{soc}(B)\rightarrow \R$ is a linear isomorphism.
	
	Next suppose that $\left(\left(A,\int_A\right),\ell\right)$ satisfies HRR$_i$, and let $\alpha\in \ker(\pi)\cap A_j$ for some $0\leq j\leq \flo{d}$.  Then $\ell\alpha=0$ and hence $\alpha\in P_{j,\ell}$.  On the other hand, by HRR$_i$, it follows that if $\alpha\neq 0$ and $0\leq j\leq i$, then $(-1)^j\int_A\ell^{d-2j}\alpha^2>0$, which leads to a contradiction if $d-2j\geq 1$.  Therefore we conclude that $\pi\colon A\rightarrow B$ is injective (and hence an isomorphism) in degrees $0\leq j\leq \min\{i,\flo{d-1}\}$, and hence $h_j(A)=h_j(B)$ for those degrees.  It remains to see that $\left(\left(B,\int_B\right),\pi(\ell)\right)$ also satisfies HRR$_i$.  First assume that $0\leq i\leq \flo{d-1}$.  Fix $0\leq j\leq i$ and fix $0\neq \alpha\in P_{j,\ell}(B)=\ker(\times\ell^{d-1-2j+1}\colon B_j\rightarrow B_{d-1-j+1})$.  Then for any lift $\hat{\alpha}\in A_j$, it follows that $0\neq \hat{\alpha}\in P_{j,\ell}(A)$, and hence by HRR$_i$ for $\left(\left(A,\int_A\right)\ell\right)$ we must have 
	$$(-1)^j\int_B\ell^{d-1-2j}\alpha^2=(-1)^j\int_A\ell^{d-2j}\alpha^2>0.$$
	This implies that $\left(\left(B,\int_B\right),\pi(\ell)\right)$ satisfies HRR$_i$.  Finally if $i=\flo{d}>\flo{d-1}$, then $\left(\left(A,\int_A\right)\ell\right)$ satisfies HRP, and since that means that it also satisfies HRR$_{i-1}$, it follows from the previous argument that $\left(\left(B,\int_B\right),\pi(\ell)\right)$ must also satisfy HRR$_{i-1}$, and hence also HRP.     
\end{proof}

	Next, we define the mixed Hodge-Riemann relations.  Given a linear form $\ell_0\in A_1$ and a sequence of linear forms $\mathcal{L}=\left(\ell_1,\ldots,\ell_{d-2j}\right)\subset \left(A_1\right)^{d-2i}$ we define the $i^{th}$ mixed primitive subspace with respect to the pair $\left(\ell_0,\mathcal{L}\right)$
	$$P_{i,\ell_0}^\mathcal{L}=\ker\left(\times \ell_0\ell_1\cdots\ell_{d-2i}\colon A_i\rightarrow A_{d-i+1}\right)$$
	and define the $i^{th}$ mixed Lefschetz form with respect to $\mathcal{L}$ 
	$$\xymatrixrowsep{.5pc}\xymatrix{A_i\times A_i\ar[r] & \R\\ (\alpha,\beta)_i^\mathcal{L}\ar@{=}[r] & (-1)^i\int_A\ell_1\cdots\ell_{d-2i}\alpha\beta\\}$$

	\begin{definition}
		\label{def:MHRR}
		The AG algebra $A$ satisfies mixed SL$_i$ on a subset $U\subset A_1$ if for every $0\leq j\leq i$, and for every sequence of linear forms $\mathcal{L}=\{\ell_1,\ldots,\ell_{d-2j}\}\subset U^{d-2j}$ the $j^{th}$ mixed Lefschetz form with respect to $\mathcal{L}$ is non-degenerate or equivalently, the $j^{th}$ mixed Lefschetz multiplication map 
		$$\times\ell_1\cdots\ell_{d-2j}\colon A_j\rightarrow A_{d-j}$$
		is an isomorphism.  
		
		The oriented AG algebra $\left(A,\int_A\right)$ satisfies the mixed HRR$_i$ on $U\subset A_1$ if for every $0\leq j\leq i$, and for every sequence of linear forms $\left(\ell_0,\mathcal{L}=\{\ell_1,\ldots,\ell_{d-2j}\}\right)\subseteq U^{d-2j+1}$, the $j^{th}$ Lefschetz form with respect to $\mathcal{L}$ is  positive definite on the $j^{th}$ primitive subspace with respect to $(\ell_0,\mathcal{L})$, i.e. 
		$$\left(\alpha,\alpha\right)_j^\mathcal{L}>0, \ \begin{array}{l} \forall 0\leq j\leq i\\ \forall 0\neq \alpha\in P_{j,\ell_0}^{\mathcal{L}}\\ \end{array}$$
		It has HRP on $U$ if it satisfies HRR$_{\flo{d}}$ on $U$.
	\end{definition}
	
	\begin{remark}
			\label{rem:Cattani}
		Note that if $\left(A,\int_A\right)$ satisfies the mixed HRR$_i$ on $U$, then it must also satisfy the ordinary HRR$_i$ on $U$ by specializing all linear forms to a single one:  $\ell_0=\ell_1=\cdots=\ell_{d-2j}=\ell$, $0\leq j\leq i$.  According to Cattani \cite{C}, the converse should also hold in the case where $U$ is a convex cone and $i=\flo{d}$, although his proof is highly non-trivial; see \Cref{sec:Conclude} for further comments.  The following example shows that in general, mixed HRR$_i$ is stronger than ordinary HRR$_i$.
	\end{remark}

\begin{example}
	\label{ex:X3Y3}
	Define the oriented AG algebra of socle degree $d=3$
	$$A=\frac{\R[x,y]}{(x^3-y^3,xy)}, \ \int_Ax^3=1.$$
	Define for $i=0,\ldots,5$ $\ell_i=a_ix+b_iy$ for $(a_i,b_i)\in \R^2$.  Then in degree $i=0$, for the sequence $(\ell_0,\underbrace{\ell_1,\ell_2,\ell_3}_{\mathcal{L}_0})$ the $0^{th}$ mixed primitive subspace 
	$$P_{0,\ell_0}^{\mathcal{L}_0}=\ker\left(\times\ell_0\ell_1\ell_2\ell_3\colon A_0\rightarrow A_4\right)=\langle 1\rangle$$
	and plugging $\alpha=1$ into the $0^{th}$ mixed Lefschetz form for $\mathcal{L}_0$ yields
	$$(1,1)_0^{\mathcal{L}_0}=\int_A\ell_1\ell_2\ell_31^2=a_1a_2a_3+b_1b_2b_3.$$
	In degree $i=1$, for the sequence $(\ell_4,\underbrace{\ell_5}_{\mathcal{L}_1})$ the $1^{st}$ mixed primitive subspace 
	$$P_{1,\ell_4}^{\mathcal{L}_1}=\ker\left(\times\ell_4\ell_5\colon A_1\rightarrow A_3\right)=\left\langle \alpha=b_4b_5x-a_4a_5y\right\rangle$$
	and plugging $\alpha$ into the $1^{st}$ mixed Lefschetz form for $\mathcal{L}_1$ yields
	$$(\alpha,\alpha)_1^{\mathcal{L}_1}=(-1)\int_A\ell_5\alpha^2=-a_5b_5\left(b_4^2b_5+a^2_4a_5\right).$$
	If we set $a_i=a$ and $b_i=b$ for all $i=0,\ldots,5$, then the set of linear forms on which $A$ satisfies the ordinary HRR$_1$ (and hence HRP) is
	$$U=\left\{ax+by \ | \ a^3+b^3>0, \ ab<0\right\}=\underbrace{\left\{ax+by \ | \ b>-a>0\right\}}_{V_1}\sqcup\underbrace{\left\{ax+by \ | \ a>-b>0\right\}}_{V_2}$$
	which is the disjoint union of two convex polyhedral cones $V_1$ and $V_2$.  However $A$ only satisfies the mixed HRR$_1$ (or mixed HRP) on \emph{either} $V_1$ \emph{or} $V_2$, but not on their union.  	 
\end{example}
	The following is a mixed analogue of \Cref{lem:HRRSL}.

		\begin{lemma}
			\label{lem:MHRRMSL}
			If $\left(A,\int_A\right)$ satisfies mixed HRR$_i$ on $U\subset A_1$, then $A$ satisfies mixed SL$_i$ on $U$.
		\end{lemma}
	\begin{proof}
	Suppose that $A$ does not satisfy mixed SL$_i$ on $U\subset A_1$.  Then for some $0\leq j\leq i$ and for some sequence of linear forms $\mathcal{L}=(\ell_1,\ldots,\ell_{d-2j})\subset U^{d-2j}$, there exists nonzero $\alpha\in \ker(\times\ell_1\cdots\ell_{d-2j}\colon A_j\rightarrow A_{d-j})$.  Then for any additional choice of $\ell_0\in U$, we see that $0\neq \alpha\in P_{j,\ell_0}^{\mathcal{L}}$, and 
	$$\int_A\ell_1\cdots\ell_{d-2j}\alpha^2=0.$$
	Hence $\left(A,\int_A\right)$ does not satisfy HRR$_i$ on $U$.
	\end{proof}

Just as in the ordinary case, the converse of \Cref{lem:MHRRMSL} is not true; in fact since mixed HRR$_i$ specializes to ordinary HRR$_i$, it follows that the the ring in \Cref{ex:SLnotHRR}, namely $A=\R[x,y]/(x^2-y^2,xy)$, also cannot satisfy mixed HRR$_1$ on any set $U$.  On the other hand, it satisfies mixed SL$_1$ on the standard open cone $U=\left\{ax+by \ | \ (a,b)\in\R^2_{>0}\right\}$ (but not its closure, since $xy=0$ in $A$).

	\begin{lemma}
	\label{lem:MPD}
	If $A$ satisfies mixed SL$_i$ on $U$, then for each $0\leq j\leq i$ and for each sequence of linear forms $\left(\ell_0,\mathcal{L}=(\ell_1,\ldots,\ell_{d-2j})\right)\subset U^{d-2j+1}$, $\dim_\R\left(P_{j,\ell_0}^{\mathcal{L}}\right)=h_j-h_{j-1}$ and there is an orthogonal decomposition with respect to the $j^{th}$ mixed Lefschetz form
	$$A_j=P_{j,\ell_0}^{\mathcal{L}}\oplus \ell_0(A_{j-1}).$$
\end{lemma}
\begin{proof}
	Fix $0\leq j\leq i$, fix $\alpha\in A_j$, fix some  $\left(\ell_0,\mathcal{L}=(\ell_1,\ldots,\ell_{d-2j})\right)\subset U^{d-2j+1}$ and consider the element $y=\ell_0\ell_1\cdots\ell_{d-2j}\alpha\in A_{d-j+1}$.  If $y=0$, then $\alpha\in P_{j,\ell_0}^{\mathcal{L}}$.  Otherwise, we may assume that $y\neq 0$ and $j>0$ and by mixed SL$_{j-1}$, for each choice of $\ell\in U$, there exists $\beta\in A_{j-1}$ such that $$\ell_0\ell_1\cdots\ell_{d-2j}\ell\beta=y=\ell_0\ell_1\cdots\ell_{d-2j}\alpha$$
	and hence $\alpha-\ell\beta\in P_{j,\mathcal{L}}$.  If we choose $\ell=\ell_0$, then for any $\alpha\in P_{j,\ell_0}^{\mathcal{L}}$ and any $\beta\in A_{j-1}$, we have $$(\alpha,\ell_0\beta)_j^{\mathcal{L}}=\int_A\ell_1\cdots\ell_{d-2j}\alpha\ell_0\beta=\int_A\ell_0\ell_1\cdots\ell_{d-2j}\alpha\beta=0.$$
	This shows that the subspace $P_{j,\ell_0}^{\mathcal{L}}$ and $\ell_0(A_{j-1})$ are orthogonal.  Finally if $\alpha\in P_{j,\ell_0}^{\mathcal{L}}\cap \ell_0(A_{j-1})$, then mixed SL$_{j-1}$ implies that $\alpha=0$ and hence the sum is direct.  Finally, mixed SL$_i$ implies that the multiplication map $\times\ell_0\colon A_{j-1}\rightarrow A_j$ is injective for all $0\leq j\leq i$, and hence $\dim_\R\left(P_{j,\ell_0}^{\mathcal{L}}\right)=\dim_\R(A_j)-\dim_\R(\ell_0(A_{j-1}))=h_j-h_{j-1}$, as desired.
\end{proof}	

	\section{Macaulay Duality and Higher Hessians}
	\label{sec:MDHess}
	Let $R=\R[x_1,\ldots,x_n]$ and $Q=\R[X_1,\ldots,X_n]$ be polynomial rings with the standard grading $\deg(x_i)=\deg(X_i)=1$, where the lower case polynomials act on the upper case polynomials by partial differentiation, i.e. $x_i\circ F=\partial F/\partial X_i$.  Then each homogeneous polynomial $F\in Q_d$ of degree $d$ defines an ideal $\Ann(F)=\left\{f\in R \left| f\circ F=0\right.\right\}$, and we define the quotient algebra $A_F=R/\Ann(F)$.  The following was discovered by Macaulay \cite{Macaulay}; see \cite[Lemma 2.14]{IK} for a proof.
	\begin{fact}
		\label{fact:Macaulay}
For any homogeneous polynomial $F\in Q_d$ of degree $d$, the algebra $A_F=R/\Ann(F)$ is an AG algebra of socle degree $d$, with canonical orientation defined by $\int_{A_F}\colon (A_F)_d\rightarrow \R$, $\int_{A_F}\alpha=\alpha\circ F$.  Moreover, every oriented AG algebra arises this way, and two oriented AG algebras are isomorphic $A_F\cong A_G$ if and only if they are related by a linear change of coordinates, i.e. $G=\sigma \cdot F$ for some $\sigma\in \GL(n,\R)$.
	\end{fact}
From now on, we shall write $A_F$ for the oriented AG algebra $\left(A=R/\Ann(F),\int_A\right)$ where $\int_A\alpha=\alpha\circ F$ is the canonical orientation on $A$.	
	\subsection{Higher Hessians and the HRR}
	Fix an oriented AG algebra $A=A_F$ of socle degree $d$. 
	\begin{definition}
	 Fix a vector space basis $\mathcal{E}=\Bigl\{e^i_p \ | \ 1\leq p\leq h_i, \ 0\leq i\leq \flo{d}\Bigr\}$ for $A_{\leq \flo{d}}$.  Then for each $0\leq i\leq \flo{d}$ define the $i^{th}$ Hessian of $F$ with respect to $\mathcal{E}$ by 
	$$\Hess_i(F,\mathcal{E})=\left(e^i_pe^i_q\circ F\right)_{1\leq p,q\leq h_i}.$$
\end{definition}
The Hessian $\Hess_i(F,\mathcal{E})$ is a symmetric matrix with polynomial entries in $Q$, and in particular for each vector $C=(C_1,\ldots,C_n)\in\R^n$, the Hessian evaluated at $C$, denoted by $\Hess_i(F,\mathcal{E})|_C$, is a real symmetric matrix. 
	 We define its signature as the number of positive eigenvalues minus the number of negative eigenvalues, i.e.
	$$\sgn\bigl(\Hess_i(F,\mathcal{E})|_C\bigr)=\#\left(\text{eigenvalues}>0\right)-\#\left(\text{eigenvalues}<0\right).$$
	
	\begin{definition}
	For a fixed linear form $\ell=\ell(C)=C_1x_1+\cdots+C_nx_n\in A_1$, let $M^{\ell(C)}_i(\mathcal{E})$ be the $h_i\times h_i$ matrix for the (signed) $i^{th}$ Lefschetz pairing, i.e. 
	$$M^{\ell(C)}_i(\mathcal{E})=\left((-1)^i\left(e_p^i,e_q^i\right)^\ell_i\right)_{1\leq p,q\leq h_i}=\left(\int_A\ell^{d-2i}e^i_pe^i_q\right)_{1\leq p,q\leq h_i}.$$
	Equivalently, $M^{\ell(C)}_i(\mathcal{E})$ is the matrix for the $i^{th}$ Lefschetz map $\times\ell^{d-2i}\colon A_i\rightarrow A_{d-i}$ with respect to the basis $\mathcal{E}_i$ and its dual basis $\mathcal{E}^*_{d-i}$ with respect to the intersection pairing on $A$.  
	\end{definition}
	
	\begin{lemma}
		\label{fact:HessLef}
		For fixed $C=(C_1,\ldots,C_n)\in\R^n$ we have equality of matrices
		$$M^{\ell(C)}_i(\mathcal{E})=d!\Hess_i(F,\mathcal{E})|_C.$$
	\end{lemma}
\begin{proof}
	The key property is the following:  For any homogeneous $d$-form $G\in Q$, and for any linear form $\ell(C)=C_1x_1+\cdots+C_nx_n\in A_1$ we have 
	$$\ell(C)^d\circ G=d!\cdot G(C_1,\ldots,C_n).$$
	To see this note that by linearity of the multiplication map it suffices to prove this in the case where $G$ is a monomial, say $G=X_1^{a_1}\cdots X_n^{a_n}$.  In this case note that 
	$$\ell^d=(C_1x_1+\cdots+C_nx_n)^d=\frac{d!}{a_1!\cdots a_n!}C_1^{a_1}\cdots C_n^{a_n}x_1^{a_1}\cdots x_n^{a_n}+\left(\text{other monomial terms of degree $d$}\right)$$
	where the ``other monomial terms'' all annihilate $G$.  Therefore we have 
	$$\ell(C)^d\circ G=\frac{d!}{a_1!\cdots a_n!}C_1^{a_1}\cdots C_n^{a_n}x_1^{a_1}\cdots x_n^{a_n}\circ X_1^{a_1}\cdots X_n^{a_n}=d!C_1^{a_1}\cdots C_n^{a_n}=d!G(C)$$
	as claimed.
	Then to complete the proof, we simply observe that for $\ell=\ell(C)=C_1x_1+\cdots+C_nx_n$, we have  
	$$(-1)^i(e^i_p,e^i_q)_i^\ell=\int_A\ell^{d-2i}e^i_pe^i_q=\ell^{d-2i}e^i_pe^i_q\circ F=d!\cdot e^i_pe^i_q\circ F|_{C},$$
	and the result follows.
\end{proof}
It follows from \Cref{fact:HessLef} that the signature of $\Hess_i(F,\mathcal{E})|_C$ is independent of our choice of basis $\mathcal{E}$; indeed changing the basis amounts to changing  $M^{\ell(C)}_i(\mathcal{E})$ to $M^T M^{\ell(C)}_i(\mathcal{E}) M$ where $M\in \GL_n(\R)$, which does not affect the signature.
	\begin{lemma}
		\label{fact:HessHRR}
		Let $A=A_F$ be an oriented AG algebra with Hilbert function $H(A)=(h_0,h_1,\ldots,h_d)$, and let $\ell(C)=C_1x_1+\cdots+C_nx_n\in A_1$ be a linear form parametrized by some $C=(C_1,\ldots,C_n)\in\R^n$.  Then the pair $\bigl(A_F,\ell(C)\bigr)$ satisfies HRR$_i$ if and only if for each $0\leq j\leq i$ we have 
		\begin{enumerate}
			\item $\det\left(\Hess_j(F,\mathcal{E})|_C\right)\neq 0$, and
			\item $\sgn\left(\Hess_j(F,\mathcal{E})|_C\right)=\sum_{k=0}^j(-1)^k(h_k-h_{k-1})$.
		\end{enumerate}
	\end{lemma}
	\begin{proof}
		By \Cref{fact:HessLef}, we have $\Hess_j(F,\mathcal{E})|_C=\frac{1}{d!}M^{\ell(C)}_j(\mathcal{E})$, and if we choose our basis $\mathcal{E}$ according to the decomposition in \Cref{lem:primitiveD}, then the matrix $M^{\ell(C)}_j(\mathcal{E})=M^\ell_j$ has the block diagonal form
		$$M^{\ell}_j=\left(\begin{array}{c|c} M^{\ell}_j|_{P_{j,\ell}} & 0\\
			\hline 
			0 & M^{\ell}_j|_{\ell(A_{j-1})}\\ \end{array}\right).$$
		Since $M^{\ell}_j|_{\ell(A_{j-1})}=M^{\ell}_{j-1}$, we obtain the formula 
		$$\sgn\left(M^{\ell}_j\right)=\sgn\left(M^\ell_j|_{P_{j,\ell}}\right)+\sgn\left(M^{\ell}_{j-1}\right).$$
		In particular it follows from \Cref{lem:HRRSL} and \Cref{lem:primitiveD} that if $(A_F,\ell(C))$ satisfies HRR$_i$ then for each $0\leq j\leq i$, $$\sgn\left(M^\ell_j|_{P_{j,\ell}}\right)=(-1)^j\dim_\R(P_{j,\ell})=(-1)^j(h_j-h_{j-1}).$$
		Using induction on $0\leq j\leq i$, it follows that if $(A_F,\ell(C))$ satisfies HRR$_i$ then conditions (1) and (2) hold.  
		
		Conversely, assume that conditions (1) and (2) hold.  Then by (1), $(A=A_F,\ell=\ell(C))$ satisfies SL$_i$, hence the decomposition from \Cref{lem:primitiveD} and also the block matrix decomposition above still holds.  It follows from the block decomposition that we have 
		$$\sgn\left(M^{\ell}_j|_{P_{j,\ell}}\right)=(-1)^j(h_j-h_{j-1}).$$
		By the primitive decomposition in \Cref{lem:primitiveD} it follows that $\dim_\R(P_{j,\ell})=h_j-h_{j-1}$ and hence the $j^{th}$-Lefschetz form is $(-1)^j$-definite on $P_{j,\ell}$, hence $(A,\ell)$ satisfies HRR$_i$ as desired.
	\end{proof}
	
	\Cref{fact:HessLef} and \Cref{fact:HessHRR} have their mixed versions as well.  To describe these mixed versions we introduce a polarization operator.  For a fixed integer $m$, define the $\R$-linear map 
	$$\xymatrixrowsep{.5pc}\xymatrix{\operatorname{Pol}_m\colon Q_m\ar[r] & \R[X_{1,1},\ldots,X_{1,m},\ldots,X_{n,1},\ldots,X_{n,m}]_m\\ 
		\operatorname{Pol}_m(G)\ar@{=}[r] & \frac{\partial^m}{\partial t_1\cdots \partial t_m} G(t_1(X_{1,1},\ldots,X_{n,1})+\cdots+t_m(X_{1,m},\ldots,X_{n,m})).\\}$$
	
	Alternatively, $\operatorname{Pol}_m$ is the $\R$-linear map that takes the product of powers of variables $X_1^{a_1}\cdots X_n^{a_n}$ of degree $m$ to the product of elementary symmetric functions of those degrees 
	\begin{align*}
		\operatorname{Pol}_m\left(X_1^{a_1}\cdots X_n^{a_n}\right)= & a_1!\cdots a_n!e_{a_1}(X_{1,1},\ldots,X_{1,m})\cdots e_{a_n}(X_{n,1},\ldots,X_{n,m})\\
		= & \sum_{\substack{|K_i|=a_i\\ K_i\cap K_j=\emptyset\\ K_1\cup\cdots\cup K_n=[m]}}a_1!\cdots a_n!\prod_{k\in K_1}X_{1,k}\cdots\prod_{k\in K_n}X_{n,k}.
	\end{align*}	
	\begin{definition}
	Fix a vector space basis $\mathcal{E}=\{e^i_p \ | \ 1\leq p\leq h_i, \ 0\leq i\leq \flo{d}\}$ for $A_{\leq \flo{d}}$. Given any subset of linear forms $\mathcal{L}(\underline{C})=\left\{\ell_1(C_1),\ldots,\ell_{d-2i}(C_{d-2i})\right\}$, where $\ell_k=\ell_k(C_k)=C_{1,k}x_1+\cdots+C_{n,k}x_n$ and $\underline{C}=(C_1,\ldots,C_{d-2j})\in\left(\R^n\right)^{d-2j}$, define the $i^{th}$ mixed Lefschetz matrix
	$$M^{\mathcal{L}(\underline{C})}_{i}(\mathcal{E})=\left((-1)^i\left(e^i_p,e^i_q\right)^{\mathcal{L}(\underline{C}))}_i\right)_{1\leq p,q\leq h_i}=\left(\int_A\ell_1\cdots\ell_{d-2i}e^i_pe^i_q\right)_{1\leq p,q\leq h_i}.$$
	Equivalently, $M_i^{\ell(C)}(\mathcal{E})$ is the matrix for the $i^{th}$ mixed Lefschetz map $\times\ell_1\cdots\ell_{d-2j}\colon A_j\rightarrow A_{d-j}$ with respect to the basis $\mathcal{E}$ and its dual basis $\mathcal{E}^*$.	
\end{definition}

	\begin{definition}
	Define the $i^{th}$ mixed Hessian of $F$ with respect to $\mathcal{E}$ as
	$$\operatorname{MHess}_i(F,\mathcal{E})=\left(\operatorname{Pol}_{d-2i}\left(e^i_pe^i_q\circ F\right)\right)_{1\leq p,q\leq h_i}.$$
	\end{definition}
	The following is a mixed version of \Cref{fact:HessLef}.
	\begin{lemma}
		\label{fact:MHEssLef}
		For fixed $\underline{C}=\left(\left(C_{11},\ldots,C_{n1}\right),\ldots,\left(C_{1,d-2j},\ldots,C_{n,d-2j}\right)\right)\in \left(\R^n\right)^{d-2i}$ we have equality of matrices
		$$M^{\mathcal{L}(\underline{C})}_{i}(\mathcal{E})=\operatorname{MHess}_i(F,\mathcal{E})|_{\underline{C}}.$$
	\end{lemma}
	\begin{proof}
		Similar to the proof of \Cref{fact:HessLef}, the key is the following formula:  For any homogeneous $m$ form $G\in Q_m$ and and for any collection of linear forms $\ell_1(C_1),\ldots,\ell_m(C_m)\in R_1$ where $\ell_i(C_i)=C_{1,i}x_1+\cdots+C_{n,i}x_n$ for $i=1,\ldots,m$, we have 
		$$\ell_1(C_1)\cdots\ell_{m}(C_m)\circ G=\operatorname{Pol}_d(G)|_{\underline{C}}.$$
		As before, it suffices to show this holds if $G$ is a monomial, say $G=X_1^{a_1}\cdots X_n^{a_n}$ where $a_1+\cdots+a_n=m$.  We have
		$$\ell_1(C_1)\cdots\ell_d(C_d)=\sum_{\substack{|K_i|=a_i\\ K_i\cap K_j=\emptyset\\ K_1\cup\cdots\cup K_n=[m]}}\prod_{k\in K_1}C_{1,k}\cdots\prod_{k\in K_n}C_{n,k}x_1^{a_1}\cdots x_n^{a_n}+\left(\text{other monomial terms}\right)$$
		where the ``other monomial terms'' all annihilate $G$.  Thus we obtain
		$$\ell_1(C_1)\cdots\ell_{d}(C_d)\circ G=\sum_{\substack{|K_i|=a_i\\ K_i\cap K_j=\emptyset\\ K_1\cup\cdot\cup K_n=[m]}}\prod_{k\in K_1}C_{1,k}\cdots\prod_{k\in K_n}C_{n,k}x_1^{a_1}\cdots x_n^{a_n}\circ X_1^{a_1}\cdots X_n^{a_n}=\operatorname{Pol}_d(G)|_{\underline{C}}$$
		as desired.  To complete the proof, it remains only to observe that we have 
		$$(-1)^i\left(e^i_p,e^i_q\right)_i^{\mathcal{L}(\underline{C})}=\ell_1(C_1)\cdots\ell_{d-2i}(C_{d-2i})e^i_pe^i_q\circ F=\operatorname{Pol}_{d-2i}\left(e^i_pe^i_q\circ F\right)|_{\underline{C}}$$
		and the result follows.		
	\end{proof}

	\begin{lemma}
		\label{fact:MHEssHRR}
		The oriented AG algebra $A=A_F$ satisfies the mixed HRR$_i$ on a subset $U\subset A_1$ if and only if for every $0\leq j\leq i$ and every sequence of linear forms $\mathcal{L}(\underline{C})\in U^{d-2j}$, we have 
		\begin{enumerate}
			\item $\det\left(\operatorname{MHess}_j(F,\mathcal{E})|_{\underline{C}}\right)\neq 0$, and
			\item $\sgn\left(\operatorname{MHess}_j(F,\mathcal{E})|_{\underline{C}}\right)=\sum_{k=0}^j(-1)^k(h_k-h_{k-1})$.
		\end{enumerate} 
	\end{lemma}
\begin{proof}
	The proof is along the same lines as \Cref{fact:HessHRR}.  Assume that $A=A_F$ satisfies mixed HRR$_i$ on $U\subset A_1$.  Then fixing $0\leq j\leq i$, fixing a sequence of linear forms $\left(\ell_0,\mathcal{L}(\underline{C})\right)=\left(\ell_0,\ell_1(C_1)\ldots,\ell_{d-2j}(C_{d-2j})\right)\subset U^{d-2j+1}$ and choosing the basis $\mathcal{E}$ according to the primitive decomposition in \Cref{lem:MPD}, the matrix $M^{\mathcal{L}(\underline{C})}_j(\mathcal{E})=M^{\mathcal{L}}_j$ has block decomposition 
	$$M^{\mathcal{L}}_j=\left(\begin{array}{c|c} M^{\mathcal{L}}_j|_{P_{j,\ell_0}^{\mathcal{L}}} & 0\\ 
			\hline
			0 & M^{\mathcal{L}}_j|_{\ell_0(A_{j-1})}\\ \end{array}\right)$$
	Since $M^{\mathcal{L}}_j|_{\ell_0(A_{j-1})}=M^{\mathcal{L}'}_{j-1}$ where $\mathcal{L}'=\left(\ell_0,\ell_0,\ell_1,\ldots,\ell_{d-2j}\right)\in U^{d-2j+2}$, it follows that 
	$$\sgn\left(M^{\mathcal{L}}_j\right)=\sgn\left(M^{\mathcal{L}}_j|_{P_{j,\ell_0}^{\mathcal{L}}}\right)+\sgn\left(M^{\mathcal{L}'}_{j-1}\right).$$
	In particular it follows from \Cref{lem:MHRRMSL} and \Cref{lem:MPD} that if $A_F$ satisfies mixed HRR$_i$ on $U$, then $$\sgn\left(M^{\mathcal{L}}_j|_{P_{j,\ell_0}^{\mathcal{L}}}\right)=(-1)^j(h_j-h_{j-1}).$$
	Using induction on $0\leq j\leq i$, it follows that if $A_F$ satisfies mixed HRR$_i$ on $U$, then conditions (1) and (2) hold.
	
	Conversely assume that conditions (1) and (2) hold.  Then by \Cref{fact:MHEssLef} $A_F$ satisfies mixed SL$_i$ on $U$, and hence by \Cref{lem:MPD}, for any $(\ell_0,\mathcal{L}(\underline{C})=\mathcal{L})\in U^{d-2j+1}$, the $(\ell_0,\mathcal{L})$ mixed primitive decomposition holds, the above block decomposition holds, and it follows that 
	$$\sgn\left(M^{\mathcal{L}}_j|_{P_{j,\ell_0}^{\mathcal{L}}}\right)=(-1)^j(h_j-h_{j-1}).$$
	By \Cref{lem:MPD}, it follows that $\dim_\R(P_{j,\ell_0}^{\mathcal{L}})=h_j-h_{j-1}$, and since this holds for all $0\leq j\leq i$, it follows that $A_F$ satisfies mixed HRR$_i$ on $U$. 
\end{proof}
	
	\subsection{The Two Variable Case}
	Fix a homogeneous polynomial $F=F(X,Y)\in Q_d=\R[X,Y]_d$ and let $A_F=\R[x,y]/\Ann(F)$ be its oriented AG algebra.  Then $A_F$ is a complete intersection, and its Hilbert function satisfies 
	$$H(A_F)=(1,2,3,\ldots,s(F)-1,s(F)^r,s(F)-1,\ldots,3,2,1).$$
	The number $s(F)=\max\{h_i(A_F) \ | \ 0\leq i\leq d\}$ is called the Sperner number of $F$ (or $A_F$), and is equal to the smallest degree of a polynomial in the ideal $\Ann(F)$.
	\begin{lemma}
		\label{fact:TwoMHessHRR}
		Let $A_F$ be a codimension two oriented AG algebra.  Then the pair $(A_F,\ell(C))$ satisfies ordinary HRR$_i$ if and only if for every $0\leq j\leq \min\{i,s(F)-1\}$ we have 
		$$\sgn\left(\det\left(\Hess_i(F,\mathcal{E})|_C\right)\right)=(-1)^{\flo{j+1}}.$$ 
		
		The oriented AG algebra $A_F$ satisfies mixed HRR$_i$ on a subset $U\subset A_1$ if and only if for each $0\leq j\leq \min\{i,s(F)-1\}$ and for each sequence $\mathcal{L}(\underline{C})\in U^{d-2j}$, we have 
		$$\sgn\left(\det\left(\operatorname{MHess}_i(F,\mathcal{E})|_{\underline{C}}\right)\right)=(-1)^{\flo{j+1}}.$$
	\end{lemma}
\begin{proof}
	Since the mixed case specializes to the ordinary case, it suffices to prove the mixed one.  Assume that $A_F$ satisfies mixed HRR$_i$ on $U$, fix $0\leq j\leq i$, and fix a sequence of linear forms $\mathcal{L}(\underline{C})=\mathcal{L}\in U^{d-2j}$.  Then according to the proof of \Cref{fact:MHEssHRR}, choosing $\ell_0\in U$, and choosing the mixed primitive basis $\mathcal{E}$ with respect to $(\ell_0,\mathcal{L})$, the matrix $M^{\mathcal{L}(\underline{C})}_j(\mathcal{E})=M^{\mathcal{L}}_j$ has the block decomposition 
	$$M^{\mathcal{L}}_j=\left(\begin{array}{c|c} M^{\mathcal{L}}_j|_{P_{j,\ell_0}^{\mathcal{L}}} & 0\\ 
		\hline
		0 & M^{\mathcal{L}'}_{j-1}\\ \end{array}\right).$$
	Note that for $0\leq j\leq \min\{i,s(F)-1\}$, we have $\dim_\R(P_{j,\ell_0}^{\mathcal{L}})=h_j-h_{j-1}=1$, whereas for $j\geq s(F)$ we must have $P_{j,\ell_0}^{\mathcal{L}}=0$.  It follows that for $0\leq j\leq \min\{i,s(F)-1\}$,  
	$$\sgn\left(M^{\mathcal{L}}_j|_{P_{j,\mathcal{L}}}\right)=(-1)^j=\sgn\left(\det\left(M^{\mathcal{L}}_j|_{P_{j,\mathcal{L}}}\right)\right)=\frac{\sgn\left(\det\left(M^{\mathcal{L}}_j\right)\right)}{\sgn\left(\det\left(M^{\mathcal{L'}}_{j-1}\right)\right)}.$$
	Inductively, assuming that $\sgn\left(\det\left(M^{\mathcal{L}'}_{j-1}\right)\right)=(-1)^{\flo{j}}$, it follows that 
	$$\sgn\left(\det\left(M^{\mathcal{L}}_j\right)\right)=(-1)^{j-\flo{j}}=(-1)^{\flo{j+1}}.$$
	Conversely assume that for all $0\leq j\leq \min\{i,s(F)-1\}$, and for all $\mathcal{L}(\underline{C})=\mathcal{L}\in U^{d-2j}$, $$\sgn\left(\det\left(M^{\mathcal{L}}_j\right)\right)=(-1)^{\flo{j+1}}.$$  
	Then $A_F$ at least satisfies mixed SL$_i$ on $U$, hence for each $\ell_0\in U$, the $(\ell_0,\mathcal{L})$ primitive decomposition holds, and hence the above block decomposition holds, and we have 
	$$\sgn\left(M^{\mathcal{L}}_j|_{P_{j,\ell_0}^{\mathcal{L}}}\right)=\frac{\sgn\left(\det\left(M^{\mathcal{L}}_j\right)\right)}{\sgn\left(\det\left(M^{\mathcal{L}'}_{j-1}\right)\right)}=(-1)^{\flo{j+1}}(-1)^{\flo{j}}=(-1)^j.$$
	In the case where $i\geq j\geq s(F)$, we must have $P_{j,\ell_0}^{\mathcal{L}}=0$, and hence the Hodge-Riemann conditions are vacuously satisfied.  It follows that $A_F$ satisfies mixed HRR$_i$ on $U$.
\end{proof}
For fixed $i$, set $\mathcal{E}_i=\left\{x^py^{i-p} \ | \ 0\leq p\leq i\right\}$,  the monomial basis for $R_i$.  Note that if $i\leq s(F)-1$, then $\mathcal{E}_i$ is also a basis for $\left(A_F\right)_i$.  The following result is a straightforward computation whose proof is left to the reader.  
\begin{lemma}
	\label{lem:TwoHess}
	If ${\dsp F=\sum_{k=0}^d\binom{d}{k}c_kX^kY^{d-k}\in Q_d}$, and $0\leq i\leq s(F)-1$, then
	$$\Hess_i(F,\mathcal{E})=(d)_{2i}\sum_{m=0}^{d-2i}\binom{d-2i}{m}X^mY^{d-2i-m}\left(c_{m+q+p}\right)_{0\leq p,q\leq i}$$
	where $(d)_{2i}=d(d-1)\cdots (d-2i+1)$.  Moreover the mixed Hessian, over the mixed polynomial ring $\R[X_1,\ldots,X_{d-2i},Y_1,\ldots,Y_{d-2i}]$, is 
	$$\operatorname{MHess}_i(F,\mathcal{E})=d!\sum_{K\subset [d-2i]}\prod_{k\in K}X_k\prod_{k\notin K}Y_k\left(c_{|K|+p+q}\right)_{0\leq p,q\leq i}$$
	where the sum is over all subsets $K$ of $[d-2i]=\{1,\ldots,d-2i\}$. 
\end{lemma}
For each $0\leq i\leq \flo{d}$, let $P_i$ denote the permutation matrix for the permutation $\{p\mapsto i-p \ | \ 0\leq p\leq i\}$.  Define the $i^{th}$ permuted (mixed) Hessian to be the product $P_i\Hess_i(F,\mathcal{E})$ ($P_i\operatorname{MHess}_i(F,\mathcal{E})$).  Since $\det(P_i)=(-1)^{\flo{i+1}}$ we obtain the following positivity criterion for HRR in terms of permuted (mixed) Hessians.
\begin{lemma}
	\label{lem:PermutedHess}
	Let $A_F$ be a codimension two oriented AG algebra.  Then the pair $(A_F,\ell(C))$ satisfies HRR$_i$ if and only if   
	$$\det\left(P_j\Hess_i(F,\mathcal{E})|_C\right)>0, \ \forall 0\leq j\leq i.$$
	
	The oriented AG algebra $A_F$ satisfies mixed HRR$_i$ on $U$ if and only if 
	$$\det\left(P_j\operatorname{MHess}_j(F,\mathcal{E})|_{\underline{C}}\right)>0, \ \begin{array}{l} \forall 0\leq j\leq i\\ \forall \mathcal{L}(\underline{C})\in U^{d-2j}\\ \end{array}$$
\end{lemma}
It follows from \Cref{lem:TwoHess} that the permuted (mixed) Hessians are monomial combinations of the Toeplitz matrices $\left(c_{m+i+q-p}\right)_{0\leq p,q\leq i}$, $0\leq m\leq d-2i$, which are the focus of the next section.

\section{Higher Lorentzian Polynomials}
\label{sec:HLor}
Throughout the section we work with homogeneous polynomials $F=F(X,Y)\in \R[X,Y]$. We set 
$F=\binom{d}{0}c_0Y^d+\binom{d}{1}c_1XY^{d-1}+\cdots+\binom{d}{d-1}c_{d-1}X^{d-1}Y+\binom{d}{d}c_dX^d$
 and we denote $A_F=R/\Ann(F)$  the associated AG quotient algebra of the polynomial ring $R=\R[x,y]$ with Hilbert function $h_i(A_F)=\dim (A_F)_i$.
 
We will also consider Toeplitz matrices. A {\em Toeplitz matrix} is a matrix in which the entries on each diagonal line parallel to the main diagonal are constant. More formally, a matrix $A$ is Toeplitz if and only if $a_{ij}=c_{i-j}$ for some sequence $\{c_k\}_{k\in\Z}$.  We shall use the notation $\mathcal{T}(m,n)\subset\mathcal{M}(m,n)$ to denote subspace of real $m\times n$ Toeplitz matrices inside the Euclidean space of all $m\times n$ matrices.  In the following three subsections we introduce three families of polynomials: strictly $i$-Lorentzian, $i$-Lorentzian, and normally stable polynomials.


	\subsection{Strictly $i$-Lorentzian Polynomials}
	We introduce our higher Lorentzian polynomials below.  Following Br\"and\'en-Huh, we define their strict versions first, then define their non-strict versions as the limits of the strict ones. 
	\begin{definition}
		\label{def:Lorentzian}
		The polynomial ${\dsp F=\sum_{k=0}^d\binom{d}{k}c_kX^kY^{d-k}}$ is called strictly Lorentzian of order $i$ or \emph{strictly $i$-Lorentzian} if for each $0\leq j\leq i$, the Toeplitz determinants satisfy
		$$\det\left(c_{m+j+q-p}\right)_{0\leq p,q\leq j}>0, \ \forall 0\leq m\leq d-2j.$$
		The set of strictly $i$-Lorentzian polynomials of degree $d$ is denoted by $\Lor(i)_d\subset Q_d$.
		\end{definition}
		
	\begin{definition}	
		Define the Lorentzian polynomials of order $i$ or the \emph{$i$-Lorentzian polynomials} to be the limits of the strictly $i$-Lorentzian polynomials.  The set of $i$-Lorentzian polynomials of degree $d$ is denoted by $L(i)_d\subset Q_d$, and it is the closure of $\Lor(i)_d$ in the Euclidean space $Q_d\cong \R^{d+1}$. 
	\end{definition}
	\begin{remark}
		\label{rem:swap}
		The transformation $c_k\mapsto c_{d-k}$ does not affect the determinant conditions in \Cref{def:Lorentzian} since 
		$$\det\left(c_{d-m-j-q+p}\right)_{0\leq p,q\leq j}=\det\left(c_{(d-m-j)+q-p}\right)^t_{0\leq p,q\leq j}$$
		where $M^t$ denotes the transpose of matrix $M$, it follows that $F(X,Y)$ is (strictly) $i$-Lorentzian if and only if $F(Y,X)$ is (strictly) $i$-Lorentzian.  
	\end{remark}
\begin{definition}\label{def:phi}
	For a polynomial ${\dsp F=\sum_{k=0}^d\binom{d}{k}c_kX^kY^{d-k}}$
 and for each $0\leq i\leq \flo{d}$, define a rectangular $(i+1)\times (d-i+1)$ Toeplitz matrix 
	$$\phi^i_d(F)=\left(\begin{array}{ccc} c_i & \cdots & c_d\\ \vdots & \ddots & \vdots\\ c_0 & \cdots & c_{d-i}\\ \end{array}\right)=\left(c_{i+q-p}\right)_{\substack{0\leq p\leq i\\ 0\leq q\leq d-i\\}}.$$
	This defines a linear isomorphism $\phi^i_d\colon Q_d\rightarrow\mathcal{T}(i+1,d-i+1)$ from the space of homogeneous bivariate forms of degree $d$ to the space of $(i+1)\times(d-i+1)$ Toeplitz matrices.
\end{definition}
	
	\begin{lemma}
		\label{fact:rank}
		The rank of $\phi^i_d(F)$ is equal to $h_i(A_F)$, which is equal to the minimum of $i+1$ or $s(F)$, i.e.
		$$\operatorname{rank}(\phi^i_d(F))=h_i(A_F)=\min\{i+1,s(F)\}.$$
	\end{lemma}
	\begin{proof}
		Fix $0\leq i\leq \flo{d}$, and define the matrix 
		$$\psi^i_d(F)=\left([x^py^{i-p}\circ F]_{X^qY^{d-i-q}}\right)_{\substack{0\leq p\leq i\\ 0\leq q\leq d-i\\}}$$
		where $[G]_{X^qY^{d-i-q}}$ means the coefficient of $X^qY^{d-i-q}$ in the monomial expansion of the polynomial $G$.
		Using the shorthand notation $(m)_n=m\cdots (m-n+1)$, one can easily see that we have 
		$$\psi^i_d(F)=\left(\binom{d}{p+q}c_{p+q}\cdot (p+q)_p\cdot (d-p-q)_{i-p}\right)_{\substack{0\leq p\leq i\\ 0\leq q\leq d-i\\}}=\left(c_{p+q}\cdot \frac{d!}{q!\cdot (d-i-q)!}\right)_{\substack{0\leq p\leq i\\ 0\leq q\leq d-i\\}}.$$
		Therefore letting $D_i=\operatorname{diag}\left(\frac{d!}{q!\cdot (d-i-q)!}\right)_{0\leq q\leq d-i}$ and $P_i$ the permutation matrix for the permutation $\{p\mapsto i-p \ | \ 0\leq p\leq i\}$, we see that 
		$$\psi^i_d(F)=P_i\cdot \phi^i_d(F)\cdot D_i.$$
		It follows that $\psi^i_d(F)$ and $\phi^i_d(F)$ have the same rank, but the rank of $\psi^i_d(F)$ is equal to the dimension of the space of $d-i$ forms $\left\{x^py^{i-p}\circ F \ | \ 0\leq p\leq i\right\}$ which is equal to $h_i(A_F)$, which is equal to either $i+1$ or else $s(F)$, whichever is smaller.
	\end{proof}
	
	For an $m\times n$ matrix $M\in \mathcal{M}(m,n)$, and for any $1\leq k\leq \min\{m,n\}$, and for any $k$-subsets $I\subset [m]$ and $J\subset [n]$ let $M_{IJ}$ be the $k\times k$ submatrix obtained from $M$ by taking rows indexed by $I$ and columns indexed by $J$.  Say that a minor is \emph{consecutive} if the subsets $I$ and $J$ consist of consecutive integers.
	Recall that $M$ is called totally positive (respectively non-negative) if every minor $M_{IJ}$ has positive (respectively non-negative) determinant.  The following result was first proved by Fekete in 1913; see Ando's comprehensive survey article \cite[Theorem 2.5]{Ando}.
	\begin{fact}
		\label{fact:GP}
		The matrix $M\in\mathcal{M}(m,n)$ is totally positive if and only if every consecutive minor has positive determinant.
	\end{fact}  
	It is important to note that this result does not hold if positive is replaced with non-negative.  For example, the matrix (found in Cryer \cite{Cryer})
	$$B=\left(\begin{array}{cccc} 1 & 1 & 1 & 0\\ 1 & 1 & 1 & 1\\ 0 & 1 & 1 & 1\\ \end{array}\right)$$
	has all of its consecutive minor determinants non-negative, but the $3\times 3$ non-consecutive minor $B_{\{1,2,3\}\{1,2,4\}}$ has negative determinant. Applying \Cref{fact:GP} to the matrix $\phi^i_d(F)$ of \Cref{def:phi} we obtain:
	
	\begin{proposition}
		\label{thm:Lorphi}
		$F\in\Lor(i)_d$ if and only if $\phi^i_d(F)$ is totally positive. 
	\end{proposition}
	\begin{proof}
		Fix $0\leq j\leq i$, and let $I=\{r,r+1,\ldots,r+j\}$ and $J=\{s,s+1,\ldots,s+j\}$ be two consecutive $j+1$-subsets for some $0\leq r\leq i-j$ and some $0\leq s\leq d-i-j$.  Then we have 
		$$\left(\phi^i_d(F)\right)_{IJ}=\left(c_{i+(s-r)+q-p}\right)_{\substack{0\leq p\leq j\\ 0\leq q\leq j\\}}$$
		where $j\leq i+(s-r)\leq d-j$.  It follows that the consecutive minor determinants of $\phi^i_d(F)$ are exactly the Toeplitz determinants of \Cref{def:Lorentzian}.  Hence it follows from \Cref{fact:GP} that $F\in\Lor(i)_d$ if and only if $\phi^i_d(F)$ is totally positive.
	\end{proof}
	
	Here is another useful criterion for strictly $i$-Lorentzian related to the (ordinary) Hodge-Riemann relations.
	\begin{lemma}
		\label{thm:LorHRR}
		$F\in\Lor(i)_d$ if and only if 
		\begin{enumerate}
			\item For each $0\leq j\leq i$ and each $0\leq m\leq d-2j$ we have 
			$$j\leq s(x^m\circ F)-1$$
			\item $(A_{x^m\circ F},y)$ satisfies HRR$_i$ for all $0\leq m\leq d$.
		\end{enumerate}
	\end{lemma}
	\begin{proof}
		First, for any positive integer $r$, let $\mathcal{E}=\mathcal{E}_r=\left\{x^py^{r-p} \ | \ 0\leq p\leq r\right\}$ be the ordered monomial basis for $R_r$, and define the $r^{th}$ extended Hessian of $G$ by the formula 
		$$\widehat{\Hess}_r(G,\mathcal{E})=\left(x^{p+q}y^{2r-p-q}\circ G\right)_{0\leq p,q\leq r};$$
		Note that the $r^{th}$ extended Hessian agrees with the $r^{th}$ Hessian if and only if $r\leq s(G)-1$ if and only if the $r^{th}$ extended Hessian is non-singular at some $C\in\R^2$.
		
		Write ${\dsp F=\sum_{k=0}^d\binom{d}{k}c_kX^kY^{d-k}}$.  Then for any $0\leq m\leq d$, we have 
		$$x^m\circ F=\sum_{k=0}^{d-m}\binom{d}{k+m}c_{k+m}(k+m)_mX^kY^{d-m-k}=(d)_m\sum_{k=0}^{d-m}\binom{d-m}{k}b_kX^kY^{d-m-k},$$
		where $b_k=c_{k+m}$.  Hence we see that the $j^{th}$ extended Hessian satisfies 
		$$\widehat{\Hess}_{j}(x^m\circ F,\mathcal{E})|_{(0,1)}=\frac{d!}{(d-2j)!}\cdot \left(c_{p+q+m}\right)_{0\leq p,q\leq j}$$
		and hence the $j^{th}$ permuted extended Hessian is
		$$P_j\widehat{\Hess}_j(x^m\circ F,\mathcal{E})|_{(0,1)}=\frac{d!}{(d-2j)!}\cdot \left(c_{j+q+m-p}\right)_{0\leq p,q\leq j}.$$
		
		To prove the assertion then, we first assume that $F\in\Lor(i)_d$, which means that 
		$$\det(c_{j+q+m-p})_{0\leq p,q\leq j}>0$$ 
		for each $0\leq j\leq i$ and for each $0\leq m\leq d-2j$.  In particular we see that for each $0\leq j\leq i$ and for all $0\leq m\leq d-2j$, the $j^{th}$ permuted extended Hessian of $x^m\circ F$ is nonsingular at $(0,1)$, and hence agrees with the $j^{th}$ permuted Hessian and we have $h_j(x^m\circ F)=j+1$, and hence that $j+1\leq s(x^m\circ F)$ which is (1).  It also follows that for fixed $0\leq m\leq d$, that 
		$$\sgn\left(\det\left(P_j\Hess_j(x^m\circ F,\mathcal{E})|_{(0,1)}\right)\right)>0, \ \forall 0\leq j\leq \max\{i,s(x^m\circ F)-1\}$$ 
		which means that $(A_{x^m\circ F},y)$ satisfies HRR$_i$, which is (2).  
		
		Conversely, assume that $F$ satisfies (1) and (2).  Then by (2), we have for every fixed $0\leq m\leq d$, and for each $0\leq j\leq \min\{i,s(x^m\circ F)-1\}$, 
		$$\sgn\left(\det\left(\Hess_j(x^m\circ F,\mathcal{E})|_{(0,1)}\right)\right)=\sgn\left(\det\left(\widehat{\Hess}_j(x^m\circ F,\mathcal{E})|_{(0,1)}\right)\right)=(-1)^{\flo{j+1}}.$$
		But by (1), we have $j\leq s(x^m\circ F)-1$ for every $0\leq j\leq i$ and every $0\leq m\leq d-2j$.  Putting these two together implies that for every $0\leq j\leq i$ and every $0\leq m\leq d-2j$
		$$\det\left(\left(c_{j+m+q-p}\right)_{0\leq p,q\leq j}\right)>0$$
		which means that $F$ is strictly Lorentzian of order $i$.
	\end{proof}

	\begin{lemma}
		\label{fact:xy}
		If $F\in\Lor(i)_d$, then $x\circ F, y\circ F\in\Lor(i)_{d-1}$.
	\end{lemma}
	\begin{proof}
		Assume $F\in\Lor(i)_d$.  Then it follows from \Cref{thm:LorHRR} that $0\leq i\leq s(F)-1$.  In particular since $s(F)\leq \flo{d}$, it follows that $0\leq i\leq \flo{d-1}$.  Then for all $0\leq j\leq i$ and all $0\leq k\leq d-1-2j$, we have $j\leq s(x^k\circ (x\circ F))-1$ and for all $0\leq k\leq d-1$, we have $(A_{x^k\circ (x\circ F)},y)$ has HRR$_i$, by \Cref{thm:LorHRR}.  This means that $x\circ F$ satisfies (1) and (2) from \Cref{thm:LorHRR} and hence $x\circ F\in\Lor(i)_{d-1}$.  A similar argument (by swapping the roles of $x$ and $y$ in \Cref{thm:LorHRR}; see \Cref{rem:swap}) shows that $y\circ F\in \Lor(i)_{d-1}$.
	\end{proof}
	
	Using \Cref{thm:LorHRR} and \Cref{fact:xy}, we can derive the following equivalence.  
	\begin{theorem}
		\label{thm:LorMHRR}
		$F\in\Lor(i)_d$ if and only if $0\leq i\leq s(F)-1$ and $A_F$ satisfies mixed HRR$_i$ on the standard closed convex cone 
		$$\overline{U}=\left\{ax+by \ | \ (a,b)\in\R^2_{\geq 0}\setminus\{(0,0)\}\right\}.$$
\end{theorem}
	\begin{proof}
		By induction on degree $d\geq 0$, the base case being trivial.  For the inductive step, assume that the equivalence holds for polynomials of degree $d-1$.  Assume first that $F\in\Lor(i)_d$.  Then certainly $i\leq s(F)-1$, by \Cref{thm:LorHRR}.  We further claim that $A_F$ satisfies mixed SL$_i$ on $\overline{U}$.  Otherwise, there exists an index $0\leq j\leq i$, and there exist linear forms $\ell_1,\ldots,\ell_{d-2j}\in \overline{U}$, and a non-zero element $\alpha\in \left(A_F\right)_j$ such that $$\ell_1\cdots\ell_{d-2j}\alpha\circ F=0.$$
		Then taking $\left(\ell_0'=\ell_{d-2j},\mathcal{L}'=(\ell_1,\ldots,\ell_{d-2j-1})\right)\in\overline{U}^{d-1-2j+1}$, and setting $\alpha_x\in A_{x\circ F}$ and $\alpha_y\in A_{y\circ F}$ to be the images of $\alpha$ we see that $\alpha_x\in P_{j,\ell_0'}^{\mathcal{L}'}(x\circ F)$ and $\alpha_y\in P_{j,\ell_0'}^{\mathcal{L}'}(x\circ F)$.
		By \Cref{thm:LorHRR}, we have that $(A_F,y)$ satisfies HRR$_i$ and, reversing the roles of $x$ and $y$, also $(A_F,x)$ satisfies HRR$_i$.  It follows from \Cref{lem:Key} that $(A_F)_j=(A_{x\circ F})_j=(A_{y\circ F})_j$ and hence in particular, $\alpha_x$ and $\alpha_y$ must each be nonzero in their respective algebras.  By \Cref{fact:xy}, both $x\circ F$ and $y\circ F$ are in $\Lor(i)_d$, and hence by our inductive hypothesis, $A_{x\circ F}$ and $A_{y\circ F}$ both satisfy mixed HRR$_i$ on $\overline{U}$.  Therefore we must have that  	
		$$(-1)^j\cdot \ell_1\cdots\ell_{d-2j-1}\cdot\alpha^2\circ (x\circ F)>0, \ \text{and} \ (-1)^j\cdot \ell_1\cdots\ell_{d-2j-1}\cdot\alpha^2\circ (y\circ F)>0.$$
		On the other hand, writing $\ell_{d-2j}=ax+by$ for some $a,b\geq 0$ not both zero, we have 
		$$0=(-1)^j\cdot \ell_1\cdots\ell_{d-2j}\cdot\alpha^2\circ F=(-1)^j\left(a\cdot \ell_1\cdots\ell_{d-1-2j}\alpha^2(x\circ F)+b\cdot \ell_1\cdots\ell_{d-1-2j}\alpha^2\circ (y\circ F)\right)$$
		which is a contradiction; therefore $A_F$ must satisfy mixed SL$_i$ on $\overline{U}$.  Next we want to show that $A_F$ must satisfy mixed HRR$_i$ on $\overline{U}$.  Fix an index $0\leq j\leq i$ and let $\mathcal{E}=\left\{x^py^{j-p} \ | \ 0\leq p\leq j\right\}$ be the standard monomial basis for $A_j$, and define the polynomial 
		$$D_j(\underline{X},\underline{Y})=\det\left(P_j\operatorname{MHess}_j(F,\mathcal{E})\right)$$
		the determinant of the permuted mixed Hessian matrix.
		Since $(A_F,y)$ satisfies HRR$_i$, it follows that for all $0\leq j\leq \min\{i,s(F)-1\}$, we have $D_j(\underline{0},\underline{1})>0$.  Since $\mathcal{L}(\underline{0},\underline{1})=(y,y,\ldots,y)\in \overline{U}^{d-2j}$, and $A_F$ satisfies mixed SL$_i$ on $\overline{U}$, it follows that 
		\[
		D_j(\underline{a},\underline{b})\neq 0, \forall \mathcal{L}(\underline{a},\underline{b})\in  \overline{U}^{d-2j}.
		\]
		  Since $\overline{U}$ is a connected set, hence $\overline{U}^{d-2j}$ is connected, it follows that 
		$$D_j(\underline{a},\underline{b})>0, \ \forall \mathcal{L}(\underline{a},\underline{b})\in \overline{U}^{d-2j}.$$
		In particular this implies that $A_F$ satisfies mixed HRR$_i$ on $\overline{U}$.
		
		Conversely, assume that $i\leq s(F)-1$ and that $A_F$ satisfies mixed HRR$_i$ on $\overline{U}$.  We will show that (1) and (2) from \Cref{thm:LorHRR} hold.  From our assumption it follows that $(A_F,x)$ satisfies HRR$_i$.  By \Cref{lem:Key}, this implies that $h_j(F)=h_j(x\circ F)$ for all $0\leq j\leq \min\{i,\flo{d-1}\}$ and that $(A_{x\circ F},x)$ also satisfies HRR$_i$.  Inductively, we see that for every $0\leq j\leq i$ and for every $0\leq k\leq d-2j$, there is an equality $h_j(F)=h_j(x^k\circ F)$.  Since $h_j(G)<j+1$ if and only if $j\geq s(G)$, it follows that $j\geq s(F)-1$ if and only if $j\geq s(x^k\circ F)-1$.  Therefore since we are assuming that $0\leq j\leq i\leq s(F)-1$, it follows that $j\leq s(x^k\circ F)-1$ for all $0\leq k\leq d-2j$, which is (1) from \Cref{thm:LorHRR}.  For (2), we want to show that $(A_{x^k\circ F},y)$ has HRR$_i$ for all $0\leq k\leq d$.  Fix $0\leq k\leq d$, set $A=A_F$ and set $B=A_{x^k\circ F}$.  For each non-zero $\alpha_B\in \ker(\times y^{d-k-2j+1}\colon B_j\rightarrow B_{d-k-j+1})$, let $\alpha\in R_j$ be any homogeneous representative, and let $\alpha_A\in A_j$ be the corresponding element in $A_j$.  Then we have 
		$$0=y^{d-k-2j+1}\cdot\alpha\circ (x^k\circ F)=x^k\cdot y^{d-k-2j+1}\cdot\alpha \circ F$$
		which means that $\alpha_A\in \ker(\times x^ky^{d-k-2j+1}\colon A_j\rightarrow A_{d-j+1})$.  Therefore by mixed HRR$_i$ on $A$ we must have 
		$$(-1)^jy^{d-k-2j+1}\cdot\alpha_B^2\circ (x^k\circ F)=(-1)^j\cdot x^ky^{d-k-2j}\cdot\alpha_A^2\circ F>0.$$
		In particular we have shown that $(A_{x^k\circ F},y)$ satisfies HRR$_i$, and since $k$ was arbitrary it holds for all $k$ which is (2) from \Cref{thm:LorHRR}.  Hence it follows from \Cref{thm:LorHRR} that $F\in\Lor(i)_d$ as desired. 
	\end{proof}
	
	\begin{corollary}
		\label{fact:changeofcoord}
		$A_F$ satisfies mixed HRR$_i$ on some closed convex cone $$\overline{W}=\left\{a\ell_1+b\ell_2 \ | \ (a,b)\in\R^2_{\geq 0}\setminus\{(0,0)\}\right\}$$
		for some linearly independent linear forms $\ell_1,\ell_2\in R_1$, and $0\leq i\leq s(F)-1$ if and only if there exists $\sigma\in\GL(2,\R)$ such that $\sigma\cdot F\in \Lor(i)_d$.
	\end{corollary}
	\begin{proof}
		Suppose that $A_F$ satisfies the hypotheses.  By Macaulay duality, every linear isomorphism $\sigma^*\colon R_1\rightarrow R_1$ defines a linear change of coordinates $\sigma\colon Q_1\rightarrow Q_1$ (the transpose) that induces an isomorphism of oriented AG algebras 
		$$\sigma^*\colon \frac{R}{\Ann(\sigma\cdot F)}\rightarrow \frac{R}{\Ann(F)}.$$
		Then choosing any linear isomorphism $\sigma^*$ mapping $\overline{U}$ onto $\overline{W}$ will give the desired $\sigma\in\GL(Q_1)=\GL(2,\R)$.  In fact we can be rather precise here:  if $\ell_1=px+qy$ and $\ell_2=rx+sy$ for some real numbers $p,q,r,s\in R$ satisfying $ps-qr\neq 0$, then $\sigma^*=\left(\begin{array}{cc} p & q\\ r & s\\ \end{array}\right)$ is the linear transformation mapping $x\mapsto \ell_1$ and $y\mapsto\ell_2$.  Furthermore their Macaulay dual generators are 
		$$L_1=\frac{1}{ps-qr}\left(sX-rY\right), \ \text{and} \ L_2=\frac{1}{ps-qr}\left(-qX+pY\right)$$
		and hence $X=pL_1+rL_2$ and $Y=qL_1+sL_2$.  Then the transpose $\sigma=\left(\begin{array}{cc} p & r\\ q & s\\ \end{array}\right)$ is the linear transformation mapping $L_1\mapsto X$ and $L_2\mapsto Y$ and hence $\sigma\cdot F(X,Y)=F(pX+rY,qX+sY)=G(X,Y)$ is the Macaulay dual generator of the oriented AG algebra $A_G$ that satisfies mixed HRR$_i$ on the standard closed cone $\overline{U}$.  Since $s(G)=s(F)$, it follows that $G\in\Lor(i)_d$.  For the converse, the argument can be reversed; the details are left to the reader.
	\end{proof}

\begin{example}
	\label{ex:X3Y3Part2}
	We have seen in \Cref{ex:X3Y3} that the oriented AG algebra 
	$$A_F=\frac{\R[x,y]}{\Ann\left(F=\frac{1}{6}\left(X^3+Y^3\right)\right)}$$
	satisfies mixed HRR$_1$ on the open (nonstandard) convex cone $V_1=\left\{ax+by \ | \ b>-a>0\right\}$.
	Note that $F\notin \Lor(1)_3$ since $\phi^1_3(F)=\left(\begin{array}{ccc} 0 & 0 & 1\\ 1 & 0 & 0\\ \end{array}\right)$
	is not totally positive.
	Fix $\ell_1=-x+2y, \ \ell_2=-x+3y\in U_1$; then $A_F$ satisfies mixed HRR$_1$ on the (nonstandard) closed cone $\overline{W}=\left\{a\ell_1+b\ell_2 \ | \ (a,b)\in\R^2_{\geq 0}\setminus\{(0,0)\}\right\}\subset V_1$.  Define the linear transformation $\sigma^*=\left(\begin{array}{cc}-1 & 2\\ -1 & 3\\ \end{array}\right)$
	which maps the standard closed cone $\overline{U}=\left\{ax+by \ | \ (a,b)\in\R^2_{\geq 0}\setminus\{(0,0)\}\right\}$ isomorphically onto $\overline{W}$.  Then the transpose $\sigma=\left(\begin{array}{rr} -1 & -1\\ 2 & 3\\ \end{array}\right)$ transforms $F$ into 
	$$G(X,Y)=F(-X-Y,2X+3Y)=\binom{3}{0}26Y^3+\binom{3}{1}17XY^2+\binom{3}{2}11X^2Y+\binom{3}{3}7X^3$$
	and $G\in\Lor(1)_3$ since 
	$\phi^1_3(G)=\left(\begin{array}{ccc} 17 & 11 & 7\\ 26 & 17 & 11\\ \end{array}\right)$
	is totally positive.
\end{example}
Using \Cref{fact:changeofcoord}, we can show ordinary HRR$_i$ implies mixed HRR$_i$, although not necessarily on the same cone.

\begin{corollary}
	\label{fact:HRRLor}
	If $(A_F,\ell)$ satisfies ordinary HRR$_i$ for some $\ell\in A_1$ where $0\leq i\leq s(F)-1$, then there exists $\sigma\in \GL(2,\R)$ such that $G(X,Y)=\sigma\cdot F(X,Y)\in\Lor(i)_d$.  In particular, $A_F$ satisfies mixed HRR$_i$ on some closed cone $\overline{W}$ containing $\ell$.
\end{corollary}
\begin{proof}
	Let $A=A_F$ and let $\ell=\ell_1\in A_1$ be any linear form such that $(A,\ell)$ satisfies HRR$_i$.  Then for any other linearly independent linear form $z\in A_1$, set $\ell_2(\epsilon)=\ell_1+\epsilon z\in A_1$.  Note that since $\ell_1$ is HRR$_i$ for $A$, it follows from \Cref{lem:Key} that for each $0\leq j\leq i$ and for each $0\leq k\leq d-2j$, $j\leq s(\ell_1^k\circ F)-1$, which is (1) in \Cref{thm:LorHRR}.  Also we see that since $(A_{\ell_1^k\circ F},\ell_1)$ is HRR$_i$ (again by \Cref{lem:Key}) and by the openness of the HRR$_i$ condition, it follows that $(A_{\ell_1^k\circ F},\ell_2=\ell_1+\epsilon z)$ also satisfies HRR$_i$ for all sufficiently small $\epsilon>0$, which is (2) from \Cref{thm:LorHRR}.  Therefore, it follows from \Cref{fact:changeofcoord} that if $\sigma^*\colon R_1\rightarrow R_1$ is the linear transformation mapping $x\mapsto \ell_1$ and $y\mapsto \ell_2$, then $G(X,Y)=\sigma\cdot F(X,Y)\in\Lor(i)_d$, and the result follows.
\end{proof}

\subsection{$i$-Lorentzian Polynomials}	
	Next we turn to $i$-Lorentzian polynomials which are limits of strictly $i$-Lorentzian polynomials.  We recall that the map $\phi^i_d\colon Q_d\rightarrow \mathcal{T}(i+1,d-i+1)$ is a linear isomorphism from the space of real homogeneous bivariate polynomials of degree $d$ onto the space of real $(i+1)\times(d-i+1)$ Toeplitz matrices.  In the search for an explicit characterization of $i$-Lorentzian polynomials, our first clue is the following fact, relating totally positive (arbitrary) matrices to totally nonnegative ones, due to A. Whitney \cite{Wh}; see \cite[Theorem 2.7]{Ando}. 
	\begin{fact}
		\label{fact:Whitney}
		The closure of the subset of totally positive matrices in the Euclidean space of $m\times n$ matrices $\mathcal{M}(m,n)$ is equal to the subset of totally non-negative matrices. 
	\end{fact}
	In light of \Cref{thm:Lorphi}, we would like to prove an analogue of \Cref{fact:Whitney} for \emph{Toeplitz} matrices.  One direction is clear:  by the continuity of $\phi^i_d$ it follows directly that if $F\in L(i)_d$ then $\phi^i_d(F)$ is totally non-negative.  The converse is more complicated.  The key idea is to realize the permuted mixed Hessian as a certain finite minor of the product of two bi-infinite Toeplitz matrices, one of which is a weighted path matrix.
	
	\subsubsection{Weighted Path Matrices and Factorization of the Permuted Hessian}
	Let $\Gamma$ be a directed acyclic graph, and let $\mathcal{A}=\{A_0,A_1,\ldots\}$ and $\mathcal{B}=\{B_0,B_1,\ldots\}$ be two (possibly infinite) sets of vertices of $\Gamma$ of the same cardinality.  For any commutative ring $Q$, a $Q$-weighting of the edges of $\Gamma$ is any function $\omega\colon E_\Gamma\rightarrow Q$, and we call $\omega(e)$ the weight of edge $e\in E_\Gamma$.  For a directed path $P\colon A_p\rightarrow B_q$, define the weight of $P$ to be the product $\omega(P)=\prod_{e\in P}\omega(e)$ of weights of the edges in $P$.  The weighted path matrix (with respect to $\Gamma$, $\mathcal{A}$, $\mathcal{B}$, and $\omega$) is the matrix 
	$$W=W(\Gamma,\mathcal{A},\mathcal{B},\omega)=\left(\sum_{P\colon A_p\rightarrow B_q}\omega(P)\right)_{0\leq p,q}$$
	where the sum is over all directed paths from $A_p\in\mathcal{A}$ to $B_q\in\mathcal{B}$.
	
	Define a path system $\mathcal{P}\colon\mathcal{A}\rightarrow\mathcal{B}$ to be a collection of paths 
	$\mathcal{P}=\left\{P_i\colon A_i\rightarrow B_{\sigma(i)} \ | \ i=0,1,\ldots\right\}$.  We say that the path system is \emph{vertex disjoint} if no two paths in $\mathcal{P}$ have a common vertex.  If $\mathcal{A}$ and $\mathcal{B}$ are finite sets, or more generally, if $I\subset\mathcal{A}$ and $J\subset\mathcal{B}$ are finite subsets of the row and column indexing sets, respectively, then we may define the sign of the path system $\mathcal{P}\colon I\rightarrow J$ to be the sign of the corresponding permutation, i.e. $\sgn(\mathcal{P})=\sgn(\sigma)$.  Define the weight of the path system to be the product of paths in the system, i.e. $\omega(\mathcal{P})=\prod_{P\in\mathcal{P}}\omega(P)$.  The following result is a well known result of Lindstr\"om \cite{L} and Gessel-Viennot \cite{GV}:
	\begin{fact}
		\label{fact:LGV}
		The determinant of the weighted path matrix is
		$$\det(W)_{IJ}=\sum_{\substack{\mathcal{P}\colon I\rightarrow J\\ \mathcal{P} \ \text{vertex disjoint}\\}}\sgn(\mathcal{P})\omega(\mathcal{P}).$$
	\end{fact}
	
As an application of \Cref{fact:LGV}, take $\Gamma=(V_\Gamma,E_\Gamma)$ to be the directed graph whose vertices $V_\Gamma$ are the $\Z^2$ lattice points in the plane, and whose directed edges $E_\Gamma$ are of the form $e\colon (a,b)\shortrightarrow (a,b+1)$ (N-step) or $e\colon (a,b)\shortrightarrow (a+1,b)$ (E-step), and hence directed paths are NE lattice paths.  Let $\mathcal{A}=\left\{A_p=(-p,p) \ | \ p=0,1,\ldots\right\}$ be the lattice points along the line $y=-x$ in the second quadrant, and for a fixed positive integer $s>0$ set $\mathcal{B}=\left\{B_q=(-q,q+s) \ | \ q=0, 1,\ldots\right\}$, the lattice points along the line $y=-x+s$ in the second quadrant.  Define the mixed weighting $\omega\colon E_\Gamma\rightarrow Q=\R[X_1,\ldots,X_s,Y_1,\ldots,Y_s]$ by $\omega(e\colon (a,b)\shortrightarrow (a+1,b))=X_{a+b}$ and $\omega(e\colon (a,b)\shortrightarrow (a,b+1))=Y_{a+b}$. Then for any path $P\colon A_p\rightarrow B_q$ we have $\omega(P)=X_{i_1}\cdots X_{i_{p-q}}Y_{j_1}\cdots Y_{j_{s-{p-q}}}$, where $P$ takes E-steps at points $(a,b)$ where $a+b\in\{i_1,\ldots,i_{p-q}\}$ and N-steps where $a+b\in\{j_1,\ldots,j_{s-(p-q)}\}$.  Therefore the weighted path matrix is 
	$$W_s(\underline{X},\underline{Y})=W_s(X_1,\ldots,X_s,Y_1,\ldots,Y_s)=\left(\sum_{K\in \binom{[s]}{p-q}}\prod_{k\in K}X_k\prod_{j\notin K}Y_j\right)_{0\leq p,q},$$ 
	a bi-infinite lower triangular Toeplitz matrix.
	
\begin{lemma}
	\label{fact:Positive1}
	For every positive integer $s$, and for each choice of nonnegative real numbers\\ 
	$a_1,\ldots,a_s,b_1,\ldots,b_s\geq 0$ the real bi-infinite Toeplitz matrix 
	$$W_s(\underline{a},\underline{b})=W_s(a_1,\ldots,a_s,b_1,\ldots,b_s)=\left(\sum_{K\in \binom{[s]}{p-q}}\prod_{k\in K}a_k\prod_{j\notin K}b_j\right)_{0\leq p,q}$$
	is totally non-negative.  More precisely, 
	\begin{enumerate}
		\item if $a_1,\ldots,b_s>0$, and $J=\{j_0<j_1<\cdots<j_t\}$ is any column set, then $\det\left(W_s(\underline{a},\underline{b})_{IJ}\right)\geq 0$ with strict inequality if and only if the row set satisfies $I\subset \{j_0,\ldots,j_t+s\}$.
		\item if $a_1,\ldots,b_s\geq 0$ where $(a_i,b_i)\neq (0,0)$ for all $i$, then for every column set $J=\{j_0,\ldots,j_s\}$, there exists some row set $I\subset \{j_0,\ldots,j_s+s\}$ satisfying $\det(W_s(\underline{a},\underline{b}))_{IJ}>0$.
	\end{enumerate}  
\end{lemma} 
\begin{proof}
	The key observation here is that the vertex disjoint path systems in the weighted acyclic graph $\Gamma=\mathcal{A}\cup\mathcal{B}$ in this case must all correspond to the identity permutation, and hence in particular have $\sgn(\mathcal{P})=1$.  Thus it follows from \Cref{fact:LGV} that for every $(t+1)$-subsets $I,J\in\binom{\N}{t+1}$ of nonnegative integers, the $(t+1)\times (t+1)$ minor 
	$\det\left(W_s(\underline{X},\underline{Y})_{IJ}\right)$ is a sum of monomials with non-negative coefficients.  In fact, if $J=\{j_0<\cdots<j_t\}$ is any column set, and $I=\{i_0<\cdots<i_t\}$ is any row set satisfying $i_0\geq j_0$ and $i_t\leq j_t+s$ then there must be at least one vertex disjoint path system $\mathcal{P}\colon \mathcal{A}_I\rightarrow\mathcal{B}_J$.  Conversely if $i_0<j_0$ or $i_t>j_t+s$, then there can be no such vertex disjoint path systems, and hence $\det(W_s(\underline{X},\underline{Y})_{IJ})\equiv 0$.  
	\begin{figure}
		\includegraphics[width=7.5cm]{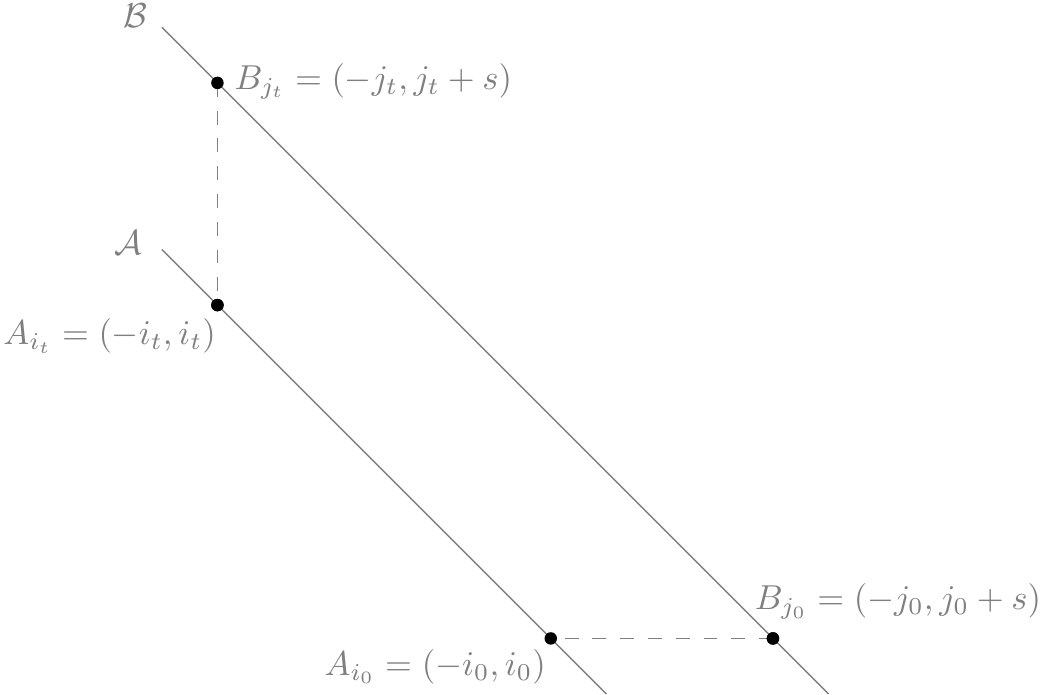}
		\caption{NE lattice paths from $\mathcal{A}$ to $\mathcal{B}$}
		\label{fig:NELP}
	\end{figure}
This implies (1).  To see (2), fix a sequence of pairs of real numbers $\left((a_1,b_1),\ldots,(a_s,b_s)\right)\in \R^2_{\geq 0}\setminus \{(0,0)\}$, and choose a sequence $\epsilon=(\epsilon_1,\ldots,\epsilon_s)$ where $0\neq \epsilon_i\in\{a_i,b_i\}$.  Then, for each $0\leq r\leq t$, starting at $B_{j_r}\in\mathcal{B}$, the sequence $\epsilon$ defines a unique path back to some $A_{i_r}\in\mathcal{A}$ as follows:  for each $0\leq k\leq s$, at a point on the line $x+y=k$, $\epsilon_k$ indicates either a E-step (if $\epsilon_k=a_k$) or an N-step (if $\epsilon_k=b_k$) taken from a unique point on the line $x+y=k-1$.  It is clear that the above path system is vertex disjoint, and we must have $I=\{i_0<\cdots<i_t\}\subset \{j_0<\cdots<j_t+s\}$, which gives (2).   
\end{proof}

	\begin{remark}
		\label{rem:FZ}
		Perhaps one should not be too surprised at the appearance of weighted path matrices here, since, according to Fomin-Zelevinsky \cite{FZ}, every totally nonnegative matrix can be realized as a weighted path matrix of some directed acyclic planar graph.
\end{remark}
The next lemma is the key to our characterization of $i$-Lorentzian polynomials; it says that the permuted mixed Hessian is a finite minor of the product of two bi-infinite Toeplitz matrices.	
	\begin{lemma}
		\label{thm:HessMHess}
		Fix ${\dsp F=\sum_{k=0}^d\binom{d}{k}c_kX^kY^{d-k}}$, fix $0\leq i\leq s(F)-1$, and let $P_i$ be the permutation matrix for the permutation $\{p\mapsto i-p \ | \ 0\leq p\leq i\}$.  Then the $i^{th}$ permuted mixed Hessian matrix is a minor of the product of bi-infinite Toeplitz matrices, specifically:
		$$P_i\cdot \operatorname{MHess}_i(F,\mathcal{E})=d!\cdot \left(\phi_d(F)\cdot W_{d-2i}(\underline{X},\underline{Y})\right)_{IJ}$$
		where $\phi_d(F)=(c_{q-p})_{0\leq p,q}$, $W_{d-2i}(\underline{X},\underline{Y})$ is as above, $I=\{0,\ldots,i\}$ and $J=\{i,\ldots,2i\}$.
	\end{lemma}
	\begin{proof}
	Computing the product of bi-infinite Toeplitz matrices gives
	\begin{align*}
		\phi_d(F)\cdot W_{d-2i}(\underline{X},\underline{Y})= & \left(c_{q-p}\right)_{0\leq p,q}\cdot\left(\sum_{K\in\binom{[d-2i]}{p-q}}\prod_{k\in K}X_k\prod_{j\notin K}Y_j\right)_{0\leq p,q}\\
		= & \left(\sum_{r\geq 0}\sum_{K\in \binom{[d-2i]}{r-q}}c_{r-p}\prod_{k\in K}X_k\prod_{j\notin K}Y_j\right)_{0\leq p,q}
	\end{align*}
and in particular for $I,J$ as above, we have 
	\begin{align*}
	\left(\phi_d(F)\cdot W_{d-2i}(\underline{X},\underline{Y})\right)_{IJ}= & \left(\sum_{r\geq 0}\sum_{K\in \binom{[d-2i]}{r-q}}c_{r-p}\prod_{k\in K}X_k\prod_{j\notin K}Y_j\right)_{\substack{0\leq p\leq i\\ i\leq q\leq 2i\\}}\\
	= & \left(\sum_{m=0}^{d-2i}\sum_{K\in\binom{[d-2i]}{m}}c_{m+q-p}\prod_{k\in K}X_k\prod_{j\notin K}Y_j\right)_{\substack{0\leq p\leq i\\ i\leq q\leq 2i\\}}\\
	= & \sum_{K\subset [d-2i]}\prod_{k\in K}X_k\prod_{j\notin K}Y_j\left(c_{|K|+i+q-p}\right)_{0\leq p,q\leq i}\\
	= & P_i\cdot \operatorname{MHess}_i(F,\mathcal{E}),
\end{align*}
as desired.
\end{proof}

\begin{remark}
	\label{rem:biinf}
	Questions of convergence naturally arise when multiplying two general bi-infinite matrices, as the entries of the product can be infinite series.  In the case above, however, both bi-infinite matrices $\phi_d(F)$ and $W_{d-2i}(\underline{X},\underline{Y})$ have only finitely many nonzero entries in every row and column, and hence all convergence issues disappear.
\end{remark}
	
\subsubsection{Characterization of Higher Lorentzian Polynomials}	
	We are now in a position to characterize $i$-Lorentzian polynomials.  
	\begin{theorem}
		\label{thm:Lorentzian}
		The following are equivalent:
		\begin{enumerate}
			\item $F\in L(i)_d$,
			\item $\phi_d^i(F)$ is totally non-negative.
			\item $A_F$ satisfies mixed HRR$_i$ on the standard open convex cone
			$$U=\left\{ax+by \ | \ (a,b)\in\R^2_{>0}\right\}.$$
		\end{enumerate}
\end{theorem}
	\begin{proof}
		(1) $\Rightarrow$ (2).  Assume that $F\in L(i)_d$.  Then there exists a sequence of polynomials $F_n\in\Lor(i)_d$ such that $\lim_{n\to\infty}F_n=F$.  By \Cref{thm:Lorphi}, we know that $\phi^d_i(F_n)$ is totally positive, hence it follows from the continuity of $\phi_d^i$ that $\phi^i_d(F)$ is totally  non-negative.
		
		(2) $\Rightarrow$ (3).  Assume that $\phi_d^i(F)$ is totally non-negative.  Fix $0\leq j\leq \min\{i,s(F)-1\}$.  By \Cref{thm:HessMHess} the $j^{th}$ mixed Hessian satisfies
		$$P_j\cdot \operatorname{MHess}_i(F,\mathcal{E})=\left(\phi_d(F)\cdot W_{d-2j}(\underline{X},\underline{Y})\right)_{IJ}$$
		where $I=\{0,\ldots,j\}$ and $J=\{j,\ldots,2j\}$.  By the Cauchy-Binet theorem we have 
		$$\det\left(\left(\phi_d(F)\cdot W_{d-2j}(\underline{X},\underline{Y})\right)_{IJ}\right)=\sum_{K\subset\binom{\N}{j+1}}\det(\phi_d(F))_{IK}\cdot\det\left(W_{d-2j}(\underline{X},\underline{Y})\right)_{KJ}$$
		where the sum is over all $j+1$-subsets of $\N=\{0,1,2,\ldots\}$.  By \Cref{fact:Positive1}, we know that for $a_1,\ldots,b_{d-2j}>0$, the $KJ$ minor of the evaluated weighted path matrix satisfies $\det(W_{d-2j}(\underline{a},\underline{b}))_{KJ}\geq 0$ with strict inequality if and only if $K\subset \{j,\ldots,d-j\}$.  On the other hand, the minors of $\phi_d(F)_{IK}$ for such $K$ are precisely the consecutive final minors of $\phi_d^i(F)$ of size $j+1$.  It follows that for any $C_1=(a_1,b_1),\ldots,C_{d-2j}=(a_{d-2j},b_{d-2j})\in \R^2_{>0}$ with $\underline{C}=(C_1,\ldots,C_n)\in\left(\R^2_{>0}\right)^{d-2j}$, we have $\det\left(P_j\cdot\operatorname{MHess}_i(F,\mathcal{E})|_{\underline{C}}\right)>0$, and (3) follows by \Cref{lem:PermutedHess}.  
		
		(3) $\Rightarrow$ (1).  Assume that $A_F$ satisfies mixed HRR$_i$ on $U$.  We show that $F\in L(i)_d$ by downward induction on the rank of the matrix $\phi^i_d(F)$.  For the base case, assume that $\operatorname{rk}(\phi^i_d(F))=i+1$ so that $i\leq s(F)-1$.  Then note that $\ell_1=x+ty$ and $\ell_2=tx+y$ are linearly independent and in $U$ for all $0<t<1$.  It follows from \Cref{fact:changeofcoord} that $G_t(X,Y)=F(X+tY,tX+Y)$ is strictly Lorentzian of order $\min\{i,s(F)-1\}=i$.  Since $\lim_{t\to 0}G_t=F$, it follows that $F\in L(i)_d$ which completes the base case.
		
		For the inductive step, assume the result holds for all homogeneous forms $H$ with $\operatorname{rk}(\phi^i_d(H))>\operatorname{rk}(\phi^i_d(F))\coloneqq s$, where $s<i+1$ (note that $s=s(F)$, the Sperner number of $F$).  Then by perturbing $G_t$ above with another parameter $u$, we get a two-parameter family of polynomials 
		$$H_{t,u}=G_t+(-1)^suY^d$$
		where $s=s(F)$.  Then for all $0\leq j\leq\flo{d}$, we have 
		$$\phi^j_d(H_{t,u})=\phi^j_d(G_t)+(-1)^suE_{0,d-2j}$$
		where $E_{0,d-j}$ is the $(j+1)\times (d-j+1)$ elementary matrix with $1$ in the $(0,d-j)$ entry (upper right corner) and zeros elsewhere.  It follows that we have, for any $1\leq k\leq j+1$, and for any $k$-subsets $A\subset\{0,\ldots,j\}$ and $B\subset\{0,\ldots,d-j\}$  
		$$\det\left(\phi^j_d(H_{t,u})\right)_{AB}=\begin{cases} \det\left(\phi^j_d(G_t)\right)_{AB} & \text{if} \ 0\notin A \ \text{or} \ d-j\notin B\\
			\det\left(\phi^j_d(G_t)\right)_{AB}+(-1)^{s+k-1}u\det\left(\phi^j(G_t)\right)_{A\setminus 0, B\setminus d-j} & \text{if} \ 0\in A \ \text{and} \ d-j\in B\\ \end{cases}$$
		 In particular it is clear that for fixed $0<t<1$, and for all $u>0$ sufficiently small, $\phi^j_d(H_{t,u})$ is totally positive for $0\leq j\leq s-1$ and for $j=s$, the maximal minors of $\phi^s_d(H_{t,u})$ are either zero, corresponding to the maximal minors of $\phi^s(G_t)$, or else positive, being a $u$-multiple of a maximal minor of $\phi^{s-1}_d(G_t)$.  Since the minors of the $(j+1)\times (j+1)$ matrix $\phi^j_d(H_{t,u})$ are also minors of the matrix $\phi^i_d(H_{t,u})$ for all $0\leq j\leq i$, it follows that $\operatorname{rk}(\phi^i_d(H_{t,u}))=s+1$, and hence our induction hypothesis applies.  It remains to show that the oriented AG algebra $A_{H_{t,u}}$ satisfies mixed HRR$_i$ on $U$.  Fix $0\leq j\leq \min\{i,s(H_{t,u})-1\}=s$.  Then applying the Cauchy-Binet formula to the product minor in \Cref{thm:HessMHess}, we obtain 
		 \begin{align*}
		 	\det\left(P_j\operatorname{MHess}_j(H_{t,u},\mathcal{E})|_{\underline{C}}\right)= &  \sum_{K\in  \binom{\{j,\ldots,d\}}{j+1}} \det\left(\phi_d(H_{t,u})\right)_{IK}\det(W_{d-2j}(\underline{a},\underline{b}))_{KJ}\\
		 	= & \sum_{K\in\binom{\{j,\ldots,d\}}{j+1}}\det\left(\phi^j_d(H_{t,u})\right)_{IK}\det(W_{d-2j}(\underline{a},\underline{b}))_{KJ}
		 	\end{align*}
			where $I=\{0,\ldots,j\}$, $J=\{j,\ldots,2j\}$, and the sum is over all $(j+1)$-subsets $K\subset \{j,\ldots,d\}$; in particular the determinant of the $j^{th}$ permuted mixed Hessian is a linear combination of all maximal minor determinants of $\phi^j_d(H_{t,u})$, which we have already showed are all non-negative with at least one being positive.  It follows from \Cref{fact:Positive1} that for each    $\underline{C}=\left((a_1,b_1),\ldots,(a_{d-2j},b_{d-2j})\right)\in\left(\R^2_{>0}\right)^{d-2j}$, we have $\det\left(W_{d-2j}(\underline{a},\underline{b})\right)>0$.  It follows that $\det\left(P_j\operatorname{MHess}_j(F,\mathcal{E})|_{\underline{C}}\right)>0$, and since this holds for each $0\leq j\leq s$, this implies that $A_{H_{t,u}}$ satisfies mixed HRR$_i$ on $U$.  Therefore by our induction hypothesis, $H_{t,u}$ is $i$-Lorentzian.  Since $\lim_{(t,u)\to 0}H_{t,u}=F$ it follows that $F$ is $i$-Lorentzian as well, which is (1).       
	\end{proof}
	
	Since $\phi^i_d\colon Q_d\rightarrow\mathcal{T}(i+1,d-i+1)$ is a linear isomorphism, we obtain the following, which is an anologue of Whitney's theorem \Cref{fact:Whitney} for Toeplitz matrices, and is \Cref{introthm:Toep} from the Introduction:
	\begin{corollary}
		\label{cor:Toep}
		The closure of the set of totally positive Toeplitz matrices is equal to the set of totally nonnegative Toeplitz matrices.
	\end{corollary}
	\begin{proof}
	Let $M$ be an $m\times n$ totally nonnegative Toeplitz matrix.  We may assume, by transposing if necessary, that $m\leq n$.  Set $i=m-1$ and $d=m+n-1$.  Then by \Cref{thm:Lorentzian}, there exists $F\in L(i)_d$ such that $M=\phi^i_d(F)$.  Then by definition of $i$-Lorentzian, there exists a sequence of strictly Lorentzian polynomials $F_n\in\Lor(i)_d$ such that $\lim_{n\to\infty}F_n=F$, and by continuity of $\phi^i_d$, we also have $\lim_{n\to\infty}\phi^i_d(F_n)=\phi^i_d(F)=M$.  By \Cref{thm:Lorphi}, the matrices $M_n=\phi^i_d(F_n)$ are totally positive Toeplitz matrices, and thus we have shown that $M$ is in the closure of the set of totally positive Toeplitz matrices.  The other containment is clear since the set of totally nonnegative Toeplitz matrices is closed and hence contains the closure of the set of totally positive Toeplitz matrices.  	
	\end{proof}

	\subsection{Stable and Normally Stable Polynomials}
	\begin{definition}
		\label{def:stable}
		A homogeneous polynomial ${\dsp F=\sum_{k=0}^d\binom{d}{k}c_kX^kY^{d-k}}$ is called \emph{stable} if the univariate polynomial $f(t)=F(1,t)$ has all real non-positive roots.
		The set of all stable polynomials of degree $d$ is denoted by $S_d\subset Q_d$.
	\end{definition}
	\begin{definition}
		\label{def:normstable}
		$F$ is called \emph{normally stable} if its tilde polynomial ${\dsp \tilde{F}=\sum_{k=0}^dc_kX^kY^{d-k}}$ is stable.
		The set of normally stable polynomials of degree $d$ is denoted by $NS_d\subset Q_d$.
	\end{definition}		
	The following fact can be deduced from \cite[Theorem 2.4.1]{B}, a fact which Brenti attributes to P\'olya-Szeg\"o.
	\begin{fact}
		\label{fact:PZ}
		Every normally stable polynomial is also stable, i.e.
		$$NS_d\subset S_d.$$
	\end{fact} 
	
	According to Niculescu \cite[Remark 1.2]{N}, the following fact was written down (though not in these terms) in a manuscript by I. Newton without proof in 1707, then later proved rigorously by C. Maclaurin in 1729.
	\begin{fact}
		\label{fact:Newton}
		Every stable polynomial is $1$-Lorentzian, i.e.
		$$S_d\subset L(1)_d.$$
	\end{fact}
	\Cref{fact:Newton} was generalized to $n\geq 2$ variables by Br\"and\'en-Huh \cite{BH}; see \Cref{fact:stable}.  It turns out, however, that \Cref{fact:Newton} does not extend to $i>1$; see \Cref{ex:NSLor}.  On the other hand, we show below that \emph{normally} stable polynomials are $i$-Lorentzian for every $i\geq 0$. 
	
	The Newton-Maclaurin theorem \Cref{fact:Newton} gives necessary conditions for a polynomial to be stable in terms of its coefficients.  Two hundred fifty years later, a complete characterization of stable polynomials, in terms of its coefficients, was found by Aissen-Schoenberg-Whitney \cite{ASW} and Edrei \cite{E}; see also \cite[Theorem 2.2.4]{B}.  Out of convenience, we state their result in terms of normally stable polynomials.
	\begin{fact}
		\label{fact:ASWE}
		The homogeneous polynomial ${\dsp F=\sum_{k=0}^d\binom{d}{k}c_kX^kY^{d-k}}$ is normally stable if and only if the bi-infinite Toeplitz matrix $\phi_d(F)=\left(c_{q-p}\right)_{0\leq p,q}$
		is totally nonnegative.
	\end{fact}
	\begin{remark}
		\label{rem:PFS}
		Sequences $(c_0,\ldots,c_d)$ satisfying the hypotheses of \Cref{fact:ASWE} are called \emph{P\'olya frequency (PF) sequences} and they have been studied extensively in the combinatorics literature, e.g. \cite{B}. 
	\end{remark}

	Note that the matrix $\phi^i_d(F)$ is a finite submatrix of the bi-infinite Toeplitz matrix $\phi_d(F)$.  Hence from \Cref{fact:ASWE} together with our characterization \Cref{thm:Lorentzian}, we deduce the following result:
	\begin{corollary}
		\label{thm:NStable}
		Every normally stable polynomial is $i$-Lorentzian for all $i\geq 0$, i.e.
		$$NS_d\subset \bigcap_{i\geq 0}L(i)_d.$$
		In particular, if $F\in NS_d$ then $A_F$ satisfies mixed HRP on $U=\{ax+by \ | \ (a,b)\in\R^2_{>0}\}$.  
	\end{corollary}
	\Cref{thm:NStable} gives an easy way to construct Lorentzian polynomials and, simultaneously, oriented AG algebras satisfying HRP.
	
	\begin{example}
		\label{ex:NSLor}
		Take $\tilde{f}(t)=t(t+1)^2(t+2)=2t+5t^2+4t^3+t^4$ so that 
		$$\tilde{F}(X,Y)=Y^4+4XY^3+5X^2Y^2+2X^3Y.$$
		Then $\tilde{F}\in S_4$, hence also $\tilde{F}\in L(1)_4$, but $\tilde{F}\notin L(2)_4$ since 
		$$\phi^2_4(\tilde{F})=\frac{1}{6}\left(\begin{array}{ccc} 5 & 3 & 0\\ 6 & 5 & 3\\ 6 & 6 & 5\\ \end{array}\right)$$
		is not totally nonnegative.  On the other hand, multiplying the coefficients of $\tilde{F}$ by the binomial coefficients we obtain 
		$$F(X,Y)=\binom{4}{0}1Y^4+\binom{4}{1}4XY^3+\binom{4}{2}5X^2Y^2+\binom{4}{3}2X^3Y+\binom{4}{4}0X^4;$$
		by definition, $F\in NS_4$, and since 
		$$\phi^2_4(F)=\left(\begin{array}{ccc} 5 & 2 & 0\\ 4 & 5 & 2\\ 1 & 4 & 5\\ \end{array}\right)$$
		is  totally non-negative, it follows that $F\in L(2)_4$ as well, although $F\notin \Lor(2)_4$ because $c_4=0$.  On the other hand, since $s(F)=3$, it follows from \Cref{fact:changeofcoord} that $G_t(X,Y)=F(X+tY,tX+Y)\in \Lor(2)_4$ for all $0<t<1$. 
	\end{example}
	
	\section{Concluding Remarks and Open Questions}
	\label{sec:Conclude}
	In \cite{C}, Cattani proves that mixed HRP on a convex cone $U\subset A_1$ is equivalent to ordinary HRP on $U$, using deep results from the theory of variation of polarized Hodge structures.  It would be nice to have an elementary proof of this fact, even for $n=2$ variables.
	\begin{problem}
		\label{prob:Cattani}
		Find an elementary proof of Cattani's theorem:  In $n\geq 2$ variables, if $U\subset A_1$ a convex cone, and if $A_F$ satisfies ordinary HRP on $U$, then $A_F$ satisfies mixed HRP on $U$. 
	\end{problem}
	
	The following example, due to Chris Eur \cite{Eur} and communicated to us by Matt Larson, shows that one cannot replace HRP by HRR$_i$ for $i<\flo{d}$ in Cattani's theorem:
\begin{example}
	\label{ex:X3Y+}
	Define the oriented AG algebra of socle degree $d=4$
	$$A_{F}=\frac{\R[x,y]}{\Ann(F=X^3Y+X^2Y^2+XY^3)}.$$
	The monomial basis $\mathcal{E}=\left\{e^i_p=x^py^{i-p} \ | \ 0\leq p\leq i\right\}$ forms a basis for $A_i$ for degrees $i=0,1,2$.  For $i=0$, the $0^{th}$ ordinary Hessian is just the polynomial, i.e. 
	$$\operatorname{Hess}_0(F,\mathcal{E})=F=X^3Y+X^2Y^2+XY^3$$
	and its determinant $\det(\Hess_0(F,\mathcal{E})|_{(a,b)})=F(a,b)>0$ for all $(a,b)\in\R^2_{>0}$.  For $i=1$, the $1^{st}$ ordinary Hessian is 
	$$\Hess_1(F,\mathcal{E})=\left(\begin{array}{lr} 6XY+2Y^2 & 3X^2+4XY+3Y^2\\
		3X^2+4XY+3Y^2 & 6XY+2X^2\\ \end{array}\right)$$
	and its determinant is 
	\begin{align*}
		\det(\Hess_1(F,\mathcal{E}))= & (6XY+2Y^2)(6XY+2X^2)-(3X^2+4XY+3Y^2)^2\\
		= & -3\left((X+Y)^4+2(X^2-Y^2)^2-4X^2Y^2\right)
	\end{align*} 
	and since $(a+b)^4-4a^2b^2\geq 0$ for all $(a,b)\in\R^2$, it follows $\det(\Hess_1(F,\mathcal{E})|_{(a,b)})<0$ for all $(a,b)\in \R^2_{>0}$.  On the other hand, for $i=2$, the $2^{nd}$ ordinary Hessian is 
	$$\Hess_2(F,\mathcal{E})=\left(\begin{array}{ccc} 0 & 6 & 4\\ 6 & 4 & 6\\ 4 & 6 & 0\\ \end{array}\right)$$
	and its determinant is $\det\left(\Hess_2(F,\mathcal{E})\right)\equiv 64>0$.  It follows from \Cref{fact:TwoMHessHRR} that $A_F$ satisfies ordinary HRR$_1$ on the standard open convex cone $U=\left\{ax+by \ | \ a,b>0\right\}$, but not ordinary HRR$_2$.
		
	On the other hand we compute the Toeplitz matrices 
	$$\phi^0_4(F)=\left(\begin{array}{ccccc} 0 & \frac{1}{4} & \frac{1}{6} & \frac{1}{4} & 0\\ \end{array}\right), \ \phi^1_4(F)=\frac{1}{12}\left(\begin{array}{cccc} 3 & 2 & 3 & 0\\ 0 & 3 & 2 & 3\\ \end{array}\right), \ \phi^2_4(F)=\frac{1}{12}\left(\begin{array}{ccc} 2 & 3 & 0\\ 3 & 2 & 3\\ 0 & 3 & 2\\ \end{array}\right)$$
	the only one of which is totally nonnegative is $\phi^0_4(F)$.  Therefore it follows from \Cref{thm:Lorentzian} that $A_F$ satisfies mixed HRR$_0$ on the standard open cone $U$, but it does not satisfy mixed HRR$_1$ (nor mixed HRR$_2$) on $U$. 
\end{example}

	
Another natural problem is to try to extend the results of this paper to $n>2$ variables.  For example, \Cref{thm:Lorentzian} gives one possible definition of  higher Lorentzian polynomial in this general setting:  $F=F(X_1,\ldots,X_n)$ is \emph{$i$-Lorentzian} if the oriented AG algebra $A_F$ satisfies mixed HRR$_i$ on the standard open cone $U=\left\{a_1x_1+\cdots+a_nx_n \ | \ (a_1,\ldots,a_n)\in\R^n_{>0}\right\}$.  
	\begin{problem}
		\label{prob:2}
	Characterize $i$-Lorentzian polynomials in $n>2$ variables.	
	\end{problem}
	\appendix
	\section{A Brief Review of Lorentzian Polynomials}
	\label{sec:App}
	In this appendix, we review some of the relevant definitions and results from Br\"and\'en-Huh \cite{BH}, and others, e.g \cite{B},\cite{H},\cite{MNY}.  We start with the definition of Lorentzian polynomials \cite[Definition 2.1]{BH}.  Here unless otherwise stated $Q=\R[X_1,\ldots,X_n]$, the standard graded polynomial ring in $n$ variables, $Q_d\subset Q$ the homogeneous polynomials of degree $d$, and $P_d\subset Q_d$ the subset of polynomials with positive coefficients. 
	\begin{definition}
	\label{def:SLorBH}
	In degree $d=0$ and $d=1$, define the strictly Lorentzian polynomials $\Lor_0=P_0\cong \R_{>0}$, $\Lor_1=P_1$, and for $d=2$ define 
	$$\Lor_2=\left\{F\in P_2 \ | \ \Hess_1(F) \ \text{is nonsingular with exactly one positive eigenvalue}\right\}.$$
	For $d>2$, define the strictly Lorentzian polynomials of degree $d$ by 
	$$\Lor_d=\left\{F\in P_d \ | \ \partial_iF\in\Lor_{d-1}, \ \forall i=1,\ldots,n\right\}$$
	where $\partial_iF=x_i\circ F$, the $i^{th}$ partial derivative of $F$.
	
	The set of Lorentzian polynomials of degree $d$ is defined to be the closure of $\Lor_d$ in the Euclidean space $Q_d$, i.e. 
	$$L_d=\overline{\left(\Lor_d\right)}.$$
	\end{definition}
	
	\begin{definition}
		\label{def:BHStable}
		A polynomial $F\in Q_d$ is stable if its coefficients are nonnegative and it is either non-vanishing on $\mathcal{H}^n$ where $\mathcal{H}\subset\C$ is the open upper half plane, or else identically zero.  The set of stable polynomials is $S_d\subset Q_d$.  
	\end{definition}
	Stable polynomials are sometimes also called hyperbolic polynomials.  According to \cite{BH}, it is equivalent to say that $F$ is stable if it has nonnegative coefficients and for every $U\in\R^n_{> 0}$ and every $V\in\R^n$, the univariate polynomial $f_{U,V}(t)=F(tU-V)$ has only real roots.  The following is \cite[Proposition 2.2]{BH}:
	\begin{fact}
		\label{fact:stable}
		Every stable polynomial is Lorentzian, i.e.
		$S_d\subset L_d$.
	\end{fact}
	
	In $n=2$ variables, the following results show that \Cref{def:SLorBH} and \Cref{def:BHStable} agrees with our \Cref{def:Lorentzian} and \Cref{def:stable}; see \cite[Example 2.3]{BH}:
	\begin{fact}
		For $n=2$ variables, a homogeneous polynomial $F=\sum_{k=0}^dc_kX^kY^{d-k}$ is strictly Lorentzian if and only if its coefficients are positive and strictly ultra log concave, i.e.
		\begin{align*}
			c_i>0, & & \forall 0\leq i\leq d\\
			\left(\frac{c_i}{\binom{d}{i}}\right)^2>\left(\frac{c_{i-1}}{\binom{d}{i-1}}\right)\left(\frac{c_{i+1}}{\binom{d}{i+1}}\right) & & \forall 1\leq i\leq d-1.
		\end{align*} 
	\end{fact}

\begin{fact}
	\label{fact:BHS}
	For $n=2$ variables, a homogeneous polynomial $F=F(X,Y)\in Q_d$ is stable if and only if the univariate polynomial $f(t)=F(1,t)$ has only real nonpositive roots.
\end{fact}

A subset $J\subset\N^n$ is said to be \emph{$M$-convex} if for any $\alpha,\beta\in J$ and any index satisfying $\alpha_i>\beta_i$, there exists an index $j$ satisfying $\alpha_j<\beta_j$ and $\alpha-e_i+e_j\in J$ where $e_i$ denotes the $i^{th}$ standard coordinate vector.  For a homogeneous polynomial $F=\sum_{\alpha\in\N^n}c_\alpha X^\alpha$, define its support $\supp(F)=\left\{\alpha\in\N^n \ | \ c_\alpha\neq 0\right\}\subset \N^n$, and we say that $F$ is \emph{$M$-convex} if its coefficients are nonnegative and its support is $M$-convex; we denote by $M_d\subset Q_d$ the subset of $M$-convex polynomials of degree $d$.  Br\"and\'en has shown \cite[Theorem 3.2]{B} that every stable polynomial is $M$-convex, i.e. $S_d\subset M_d$.
The following is \cite[Theorem 2.5]{BH} and is one of the central results of that paper.
\begin{fact}
	\label{fact:LorBH}
	In degrees $d=0$, $d=1$, and $d=2$, the Lorentzian polynomials satisfy $L_d=S_d$, and for $d>2$ we have 
	$$L_d=\left\{F\in M_d \ | \ \partial_iF\in L_{d-1}\right\}.$$
\end{fact} 

In $n=2$ variables, the $M$-convexity condition is equivalent to saying that the sequence of coefficients $(c_0,\ldots,c_d)$ has \emph{no internal zeros}, meaning that whenever $c_ic_k\neq 0$ it follows that $c_j\neq 0$ for all $0\leq i<j<k\leq d$.  The following description for Lorentzian polynomials in $n=2$ variables is \cite[Example 2.26]{BH}:
\begin{fact}
	\label{fact:BHLor}
	For $n=2$ variables, a homogeneous polynomial $F=\sum_{k=0}^dc_kX^kY^{d-k}$ is Lorentzian if and only if its coefficients are nonnegative, ultra log concave, with no internal zeros, i.e.
	\begin{align*}
	c_i\geq 0, & & \forall 0\leq i\leq d\\
	\left(\frac{c_i}{\binom{d}{i}}\right)^2\geq \left(\frac{c_{i-1}}{\binom{d}{i-1}}\right)\left(\frac{c_{i+1}}{\binom{d}{i+1}}\right) & & \forall 1\leq i\leq d-1\\
	(c_0,\ldots,c_d) && \text{no internal zeros}	
\end{align*}
\end{fact}

The inequalities in \Cref{fact:BHLor} are sometimes referred to as Newton's inequalities, after his discovery that they hold for nonnegative real rooted univariate (or homogeneous stable) polynomials; see \cite{N}.  The conditions in \Cref{fact:LorBH} and \Cref{fact:BHLor} are equivalent to the total positivity, respectively total nonnegativity, of the two rowed Toeplitz matrix 
$$\phi^1_d(F)=\left(\begin{array}{cccc} \tilde{c}_1 & \tilde{c}_2 & \cdots & \tilde{c}_d\\ 
	\tilde{c}_0 & \tilde{c}_1 & \cdots & \tilde{c}_{d-1}\\ \end{array}\right), \ \ \tilde{c}_k=\frac{c_k}{\binom{d}{k}}$$
which is our condition for strictly $1$-Lorentzian, respectively $1$-Lorentzian from \Cref{thm:Lorphi} and \Cref{thm:Lorentzian}.

Related to Hodge-Riemann relations, the following is essentially \cite[Theorem 2.16]{BH}:
\begin{fact}
	If $F\in \Lor_d$ then $(A_F,\ell(C))$ satisfies HRR$_1$ for all $C\in\R^n_{>0}$.
\end{fact}
\noindent Subsequently, Murai-Nagaoka-Yazawa \cite[Theorem 3.8]{MNY} improved it with the following:
\begin{fact}
	\label{fact:SNY}
	If $F\in L_d$ then $(A_F,\ell(C))$ satisfies HRR$_1$ for all $C\in\R^n_{>0}$.
\end{fact}
Later, Huh \cite[Proposition 5]{H} proved the following result which is related to the $i=1$ case in our \Cref{thm:Lorentzian}.
\begin{fact}
	\label{fact:Huh}
	If $A_F$ satisfies mixed HRR$_i$ on $\overline{U}\cong \R^n_{\geq 0}\setminus \{0\}$ then $F$ is Lorentzian.
\end{fact}   	

\section*{Acknowledgements}
The authors are grateful to the series of annual Lefschetz Properties In Algebra, Geometry, Topology and Combinatorics workshops, some of which we each participated in, beginning with G\"ottingen (2015), and followed by meetings at Banff (2016), Mittag Leffler (2017), Levico (2018), CIRM Luminy (2019), Oberwolfach (2020), Cortona (2022), and the Fields Institute (2023).  The first author was partially supported by CIMA -- Centro de
Investiga\c{c}\~{a}o em Matem\'{a}tica e Aplica\c{c}\~{o}es,
Universidade de \'{E}vora, project UIDB/04674/2020 (Funda\c{c}\~{a}o
para a Ci\^{e}ncia e Tecnologia).  The third author was supported by NSF DMS-2101225.  The fourth author was supported by JSPS KAKENHI Grant Number JP20K03508.

\bigskip

	\bibliographystyle{alpha}

\begin{thebibliography}{999}
		\bibitem{AHK}
		Adiprasito, K., Huh, J., Katz, E.
		{\em Hodge Theory for Combinatorial Geometries}
		Ann.\ of Math.\ vol.\ 188, 381-452 (2018).
		
		\bibitem{ASW}
		Aissen, M., Schoenberg, I.J., Whitney, A.
		{\em On the generating functions of totally positive sequences I},
		Journal d’Analyse Mathematique, 2, 93--103, (1951).
		
		\bibitem{Ando}
		Ando, T.,
		{\em Totally Positive Matrices}
		Linear Algebra and its Applications, 165-219 (1987).
		
		\bibitem{Bra}
		Br\"and\'en, P., 
		{\em Polynomials with the half-plane property and matroid theory},
		Adv.\ Math.\ vol.\ 216
		(2007), 
		no.\ 1,
		302--320.
		
		\bibitem{BH}
		Br\"and\'en, P., Huh, J., 
		{\em Lorentzian polynomials},
		Ann.\ of Math.\ vol.\ 192
		(2020), 
		no.\ 3,
		821--891.
		
		\bibitem{BW}
		Br\"ande\'n, P., Wagner, D.,
		{\em A converse to the Grace–Walsh–Szegő theorem}, 
		Mathematical Proceedings of the Cambridge Philosophical Society, 147(2), 447-453, (2009).
		
		\bibitem{B}
		Brenti, F., 
		{\em Unimodal, log-concave and Pólya frequency sequences in combinatorics}, 
		Mem.\ Amer.\ Math.\ Soc., no.\ 413, 1989	
		
		\bibitem{C}
		Cattani, E.,
		{\em Mixed Lefschetz Theorems and Hodge-Riemann Bilinear Relations}, International Mathematics Research Notices, Volume 2008, (2008). 
		
		\bibitem{Cryer}
		Cryer, C.W.
		\emph{Some properties of totally positive matrices},
		Linear Algebra and its Applications, Volume 15, Issue 1, 1--25, (1976).
		
		\bibitem{E}
		Edrei, A.,
		{\em Proof of a conjecture of Schoenberg on the generating function of a totally positive sequence},
		Canad. Math. J. 5, 86-94 (1953).
		
		\bibitem{EW}
		Elias, B., Williamson, G.,
		\emph{Hodge Theory of Soergel Bimodules}
		Ann. Math. (2) 180, No. 3, 1089--1136 (2014)
		
		\bibitem{Eur}
		Eur, C.,
		\emph{Hodge-Riemann relations and Lorentzian polynomials},
		unpublished notes available at
		\begin{verbatim}
			https://people.math.harvard.edu/~ceur/notes_pdf/Eur_IntroLorentzianPolynomials.pdf
		\end{verbatim}
		
		\bibitem{FZ}
		Fomin, S., Zelevinsky, A., 
		{\em Total Positivity:  Tests and Positivity},
		Math. Intelligencer Vol. 22 Issue 1, 23-33 (2000).
		
		
		\bibitem{GV}
		Gessel, I., Viennot, G., 
		{\em Binomial determinants, paths, and hook length formulae}, 
		Adv. Math. 58 (1985), no. 3, 300–321.
		
		\bibitem{H}
		Huh, J.
		\emph{Combinatorics and Hodge Theory},
		Proceedings of the International Congress of Mathematicians, 1, (2022).
		
		\bibitem{I}
		Iarrobino, A.
		\emph{Associated Graded Algebra of a Gorenstein Artin Algebra},
		Memoirs of the AMS, Vol. 107, No. 514, Providence, RI (1994).
		
		\bibitem{IK}
		Iarrobino, A., Kanev, V., 
		{\em Power sums, Gorenstein algebras, and determinantal loci},  Appendix C by Iarrobino and Steven L. Kleiman. 
		Lecture Notes in Mathematics, 1721, Springer-Verlag, Berlin, 1999. 
		
		\bibitem{L}
		Lindstr\"om, B., 
		{\em On the vector representations of induced matroids}, 
		Bull. London Math. Soc. 5 (1973), 85–90.
		
		\bibitem{MW}
		Maeno, T., Watanabe, J.,
		{\em Lefschetz elements of artinian Gorenstein algebras and Hessians of homogeneous polynomials}, 
		Illinois J. Math.\ 53 
		(2009) 
		593--603.
		
		\bibitem{Macaulay}
		Macaulay, F.S.,
		{\em The algebraic theory of modular systems}, Cambridge Mathematical Library. Cambridge University Press, Cambridge, 1994. 
		
		\bibitem{McM}
		McMullen, P.,
		{\em On Simple Polytopes},
		Invent. Math. 113, No. 2, 419--444 (1993).
		
				
		\bibitem{MNY}
		Murai, S., Nagaoka, T., Yazawa, A.,
		{\em Strictness of the log concavity of generating polynomials of matroids},
		arXiv preprint arXiv:2003.09568v1 (2022).
		
		\bibitem{N}
		Niculescu, C.,
		{\em A new look at Newton's inequalities},
		Journal of Inequalities in Pure and Applied Mathematics, Vol. 1, Issue 2, Article 17 (2000).
		
		\bibitem{S}
		Stanley, R.P.,
		\emph{Weyl groups, the hard Lefschetz theorem, and the Sperner property},
		SIAM J. Algebraic Discrete Methods 1, 168--184 (1980).
		
		\bibitem{T}
		Timorin, V.A., 
		\emph{An Analogue of the Hodge-Riemann Bilinear Relations for Simple Convex Polytopes},
		Russ. Math. Surv. 54, No. 2, 381--426 (1999); translation from Usp. Mat. Nauk 54, No. 2, 113--162 (1999).
		
		
		
		
		\bibitem{W}
		Watanabe, J., 
		{\em A remark on the Hessian of homogeneous polynomials}, 
		in: The Curves Seminar at
		Queen’s, vol. XIII, in: Queen’s Papers in Pure and Appl. Math., vol. 119, 2000, pp. 171--178.
		
		\bibitem{Wh}
		Whitney, A. 
		{\em A reduction theorem for totally positive matrices}, 
		J. Analyse Math., 2, 88-92 (1952). 
		
		
	\end{thebibliography}

\end{document}